\documentclass[11pt,leqno]{article}
\usepackage{amsthm,amsfonts,amssymb,amsmath,oldgerm}
\usepackage{epsfig}
\numberwithin{equation}{section}

\usepackage[hmargin={1.0in,1.0in},
			vmargin={1.0in,1.0in},
			headheight=14pt,
			]{geometry}

\usepackage{appendix}
\renewcommand\appendix{\par
\setcounter{section}{0}%
\setcounter{subsection}{0}%
\setcounter{table}{0}
\setcounter{figure}{0}
\setcounter{equation}{0}
\gdef\thetable{\Alph{table}}
\gdef\thefigure{\Alph{figure}}
\gdef\theequation{\Alph{section}-\arabic{equation}}
\section*{Appendix}
\gdef\thesection{\Alph{section}}
\setcounter{section}{0}} 

\newcommand\br{\begin{remark}}
\newcommand\er{\end{remark}}
\newcommand\bp{\begin{pmatrix}}
\newcommand\ep{\end{pmatrix}}
\newcommand{\be}{\begin{equation}}
\newcommand{\ee}{\end{equation}}
\newcommand\ba{\begin{equation}\begin{aligned}}
\newcommand\ea{\end{aligned}\end{equation}}

\newcommand{\bap}{\begin{app}}
\newcommand{\eap}{\end{app}}
\newcommand{\begs}{\begin{exams}}
\newcommand{\eegs}{\end{exams}}
\newcommand{\beg}{\begin{example}}
\newcommand{\eeg}{\end{exaplem}}
\newcommand{\bpr}{\begin{proposition}}
\newcommand{\epr}{\end{proposition}}
\newcommand{\bt}{\begin{theorem}}
\newcommand{\et}{\end{theorem}}
\newcommand{\bc}{\begin{corollary}}
\newcommand{\ec}{\end{corollary}}
\newcommand{\bl}{\begin{lemma}}

\newcommand{\bd}{\begin{definition}}
\newcommand{\ed}{\end{definition}}
\newcommand{\brs}{\begin{remarks}}
\newcommand{\ers}{\end{remarks}}

\newtheorem{theorem}{Theorem}[section]
\newtheorem{proposition}[theorem]{Proposition}
\newtheorem{corollary}[theorem]{Corollary}
\newtheorem{lemma}[theorem]{Lemma}

\theoremstyle{remark}
\newtheorem{remark}[theorem]{Remark}
\theoremstyle{definition}
\newtheorem{definition}[theorem]{Definition}

\newtheorem{example}[theorem]{Example}

\newcommand{\beq}{\begin{equation}}
\newcommand{\eeq}{\end{equation}}

\title{
Modulational Instability in Equations of KdV Type
}

\author{\sc\small
Jared C. Bronski\thanks{Department of Mathematics, University of Illinois at Urbana-Champaign, Urbana, IL 61801; jared@math.uiuc.edu.}
~~~
Vera Mikyoung Hur\thanks{Department of Mathematics, University of Illinois at Urbana-Champaign, Urbana, IL 61801; vera@math.uiuc.edu.}
~~~ 
Mathew A. Johnson\thanks{Department of Mathematics, University of Kansas, 1460 Jayhawk Boulevard, 
Lawrence, KS 66045; matjohn@math.ku.edu.}
}
\begin{document}
\maketitle

%
%
%


\begin{abstract}
It is a matter of experience that nonlinear waves in dispersive media, propagating primarily in one direction, may appear periodic in small space and time scales, but their characteristics --- amplitude, phase, wave number, etc. --- slowly vary in large space and time scales. In the 1970's, Whitham developed an asymptotic (WKB) method to study the effects of small ``modulations" on nonlinear periodic wave trains. Since then, there has been a great deal of work aiming at rigorously justifying the predictions from Whitham's formal theory. We discuss recent advances in the mathematical understanding of the dynamics, in particular, the instability of slowly modulated wave trains for nonlinear dispersive equations of KdV type.
\end{abstract}

\section{Introduction}\label{s:intro}
The modulational instability, at the core, is the observation that 
nonlinear wave trains in a dispersive medium may be unstable to self modulation, 
developing nontrivial large-scale structures, which continue to evolve as they propagate.
The phenomenon is familiar to researchers in many fields. 
Moreover, it is a fundamental issue in the theory of nonlinear waves.
The history of the modulational instability of nonlinear dispersive waves is long and fascinating.
We do not have space here for a complete account of it, and we will discuss only selected high points. 
For more details, the interested reader is encouraged to read 
the excellent review article by Zakharov and Ostrovsky~\cite{Z_O_2009}.

The first understanding of the importance of the modulational instability 
arouse in hydrodynamics in the 1960's. 
Whitham~\cite{Whitham.1967} applied his newly developed theory 
of nonlinear dispersion to Stokes' periodic waves in water, 
associating the hyperbolic (or elliptic) nature of the resulting modulation equations 
with the stability (or instability, respectively) of the underlying wave;
he showed that a Stokes wave is unstable if $k h_0>1.363\dots$, 
where $k$ is the wave number of the underlying wave and $h_0$ is the undisturbed fluid depth. 
At about the same time, Benjamin and Feir~\cite{Benjamin.Feir.1967}, 
motivated by laboratory experiments~\cite{Feir.1967}, derived their eponymous instability;
they showed that a Stokes wave in deep water is unstable. 
Additionally 
they discovered the sideband nature of the instability --- 
the marginal stability of the wave number, with a band of unstable wave numbers to either side.
Corroborating results arrived nearly simultaneously, but independently,
by Lighthill~\cite{Lighthill}, Ostrovsky~\cite{Ostrovsky}, Benney and Newell~\cite{Benny.Newell}, and Zakharov \cite{Zakharov1, Zakharov2}. 

The modulational instability has been experimentally observed, either directly or indirectly,
in many physical systems.
For example, pulses in nonlinear optical fibers with so-called anomalous dispersion are   
(approximately) governed by the focusing, cubic nonlinear Schr\"odinger equation,
which exhibits the modulational instability. The importance of the effect was noted 
in early papers by Anderson and Lisak~\cite{Anderson.Lisak.1984} and Hasegawa~\cite{Hasegawa.1984}, 
where it was proposed as a method for generating short optical pulses; 
it was later observed experimentally by Tai, Hasegawa and Tomita~\cite{Tai.Hasegawa.Tomita.1986}. 
Since then, there have been a great deal of work 
observing the phenomenon in a number of optical settings; 
see \cite{Bergman.et.al.1999}, \cite{Tai.Tomita.Jewell.Hasegawa.1986}, \cite{Stegeman.et.al.1997}, \cite{Bronski.Segev.Weinstein.2001}, \cite{Segev.et.al.2000}, for instance.

The modulational instability and associated Whitham equations for integrable systems, 
such as the KdV and nonlinear Schr\"odinger equations, 
were studied by Flashka, Forest and McLaughlin~\cite{Flaschka.Forest.McLaughlin.1980} 
and by Lax and Levermore~\cite{Lax.Levermore.1983}, 
and it has grown into a large and active field of research; 
see \cite{Ercolani.Forest.McLaughlin.Montgomery.1987}, \cite{Ercolani.Forest.McLaughlin.Montgomery.1987}, \cite{Dubrovin.Novikov.1989}, \cite{Krichever.1991}, \cite{Tian.Ye.1999}, \cite{Grava.2001}, 
for instance. 
Moreover, there has been a great deal of recent work aiming at making the connection 
between the elliptic nature of the modulational equations and the instability of the underlying wave. 
To single out a few in the vein, we mention 
\cite{Collet.Eckmann.1990}, \cite{Kirrmann.Schneider.Mielke.1992}, \cite{Bridges.Rowlands.1994}, \cite{Bridges.Mielke.1995}. 
Here we will be mainly concerned with the latter program, 
how Whitham's formal theory translates into rigorous stability results. 
This has led to many applications in nonlinear dispersive equations;
see \cite{Haragus.Kapitula.2008}, \cite{Haragus.Scheel.2002}, \cite{Gallay.Haragus.2007a}, \cite{Gallay.Haragus.2007b}, \cite{Hur.Lin.2008}, \cite{Johnson.2010}, \cite{Bronski.Johnson.2010}, 
\cite{JZ}, \cite{BrJZ}, \cite{Johnson.Noble.Rodrigues.Zumbrun.2014}, \cite{Johnson.Noble.Rodrigues.Zumbrun.2013a}, \cite{Johnson.Noble.Rodrigues.Zumbrun.2013b},
for instance.

We begin in Section~2 by constructing periodic traveling waves to the generalized KdV equation. 
In Section~3 we will review Whitham's theory of modulations
and discuss implications for the modulational instability.
Section~4 contains the main analytical contributions, 
providing a rigorous mathematical program to determine the modulational instability 
of periodic traveling waves of KdV type equations. 
For illustrative purposes, in Section~5, we will discuss applications  
to the KdV, modified KdV and Schamel equations.
We will then extend the theory to equations admitting \emph{nonlocal} dispersion 
and discuss applications to the Benjamin-Ono equation, 
as well as small amplitude waves for the fractional KdV, intermediate long-wave equations
and Whitham's equation for water waves.

\section{Periodic traveling waves of generalized KdV equations}\label{s:exist}

Much of the chapter is devoted to the modulational instability of periodic traveling waves
to the generalized Korteweg-de Vries (KdV) equation
\begin{equation}\label{gkdv}
u_t=u_{xxx}+f(u)_x.
\end{equation}
Here $t\in\mathbb{R}$ is typically proportional to elapsed time 
and $x\in\mathbb{R}$ is related to the spatial variable in the predominant direction of wave propagation;
$u=u(x,t)$ is real valued, representing the wave profile or a velocity, 
and $f$ is a suitable nonlinearity. In many examples of interest, $f$ obeys a power law.
Throughout we express partial differentiation either by a subscript or using the symbol $\partial$. 

Perhaps the best known among equations of the form \eqref{gkdv} is the KdV equation
\[
u_t=u_{xxx}+(u^2)_x
\]
itself which was put forward in \cite{Boussinesq} and \cite{KdV.1895} 
to model the unidirectional propagation of surface water waves 
with small amplitudes and long wavelengths in a channel.
It has since found relevances in other situations such as Fermi-Pasta-Ulam lattices;
see \cite{FPU.1955}, for instance.
In the case of $f(u)=au^3+bu^2$, $a, b$ constants, 
\eqref{gkdv} is called the modified KdV or Gardner equation,
which models internal waves propagating in a density stratified fluid; 
see \cite{Grimshaw.2002}, for instance. 
In the case of $f(u)=|u|^{3/2}$, moreover, \eqref{gkdv} is called the Schamel equation;
see \cite{schamel.1973}, for instance.

By a traveling wave, we mean a solution which propagates at a constant speed without changing the shape.
Therefore it takes the form $u(x,t)=u(z)$, $z=x-ct$, where $c \in \mathbb{R}$ is the wave speed 
and $u$ is a $t$-independent solution of 
\[
u_t=u_{zzz}+cu_z+f(u)_z.
\]
It then reduces by quadrature to 
\begin{equation}\label{quad1}
u_{zz}+cu+f(u)=a
\end{equation}
for some constant $a\in\mathbb{R}$. 
It further reduces by quadrature to 
\begin{equation}\label{quad2}
\frac{1}{2}u_z^2=E+au-\frac{1}{2}cu^2-F(u),
\end{equation}
where $F'=f$, $F(0)=0$, for some constant $E\in\mathbb{R}$.  
We refer to \eqref{quad2} as the profile ODE corresponding to \eqref{gkdv}.
One may employ elementary phase plane analysis to infer 
the existence of periodic, homoclinic (pulses) and heteroclinic (fronts) orbits of \eqref{quad2}. 
Note that the existence and non-existence of bounded solutions of \eqref{quad2}
is determined by the effective potential
\[
V(u;a,c):=F(u)+\frac{1}{2}cu^2-au.
\]

\subsection{Some Explicit Solutions}\label{s:explicit}
In some special cases, \eqref{quad2} may be directly integrated to yield classes of explicit solutions.
In the case of $f\equiv0$ in \eqref{gkdv}, for example, namely the (linear) Airy equation, 
\eqref{gkdv} admits plane wave solutions $u(x,t)=e^{ik(x-ct)}$, 
where $k$ represents the wave number and $c$ is the wave speed,
provided that the temporal frequency $\omega:=kc$ is related to $k$ 
via the dispersion relation $\omega=-k^3$.
We will return to such linear waves in Section~3.1 when we discuss Whitham's theory of modulations. 
The modulational stability of these waves, it turns out, is trivial 
and follows directly from the dispersion relation. 
In fact, all these waves are modulationally stable. 
Therefore we will be concerned with a genuine nonlinearity,
typically $f(u)=\sigma u^{p+1}$ for some $p>0$ and $\sigma\neq 0$ real.
Thanks to the scaling invariance of \eqref{gkdv}, by the way, $\sigma$ may be taken arbitrary. 
In fact, the stability of periodic traveling waves of \eqref{gkdv} is independent of the choice of $\sigma$.

In the case of $f(u)=\frac{1}{2}u^2$, namely the KdV equation, moreover, simple phase plane analysis implies that 
periodic solutions of \eqref{quad2} exist, depending on $a$, $E$, and $c$, 
provided that the cubic polynomial $E-V(\cdot\,;a,c)$ has three distinct real roots, 
which we label as $\gamma<\beta<\alpha$.
Periodic solutions then correspond to oscillation between in the interval $[\alpha,\beta]$, 
and $a$, $E$, and $c$ are related to $\alpha$, $\beta$, and $\gamma$ via
\begin{equation}\label{trans}
E=\frac{\alpha\beta\gamma}{6}, \qquad 
a=-\frac{\alpha\beta+\beta\gamma+\alpha\gamma}{6},
\qquad c=\frac{\alpha+\beta+\gamma}{3}.
\end{equation}
In particular, \eqref{quad2} may be written as $u_z^2=(\alpha-u)(u-\beta)(u-\gamma)/3$
and $u$ may be defined implicitly via
\[
\frac{z+z_0}{\sqrt{3}}=\int_u^\alpha\frac{dw}{\sqrt{(\alpha-w)(w-\beta)(w-\gamma)}},
\]
where $z_0$ is a constant of integration associated with translational invariance 
and can be chosen so that $u(0)=\alpha$.
We then make the change of variables
\[
w=\alpha-(\alpha-\beta)\sin^2\phi, \qquad \sin\phi=\sqrt{\frac{\alpha-w}{w-\beta}},
\]
and solve the above for $u$ explicitly as
\begin{equation}\label{cnoidal1}
u(z;\alpha,\beta,\gamma,z_0)=\beta+(\alpha-\beta){\rm cn}^2\Big(\sqrt{\frac{\alpha-\gamma}{12}}~(z+z_0);m\Big),\qquad m=\frac{\alpha-\beta}{\alpha-\gamma},
\end{equation}
where ${\rm cn}(\cdot\,;m)$ denotes the Jacobi cnoidal function with elliptic modulus $m\in[0,1)$.
These cnoidal waves define a four-parameter family of periodic solutions of the KdV equation.
Note that the period is $T=4\sqrt{3}K(m)/\sqrt{\alpha-\beta}$, 
where $K$ is the complete elliptic integral of the first kind\footnote{
Note that ${\rm cn}(z;m)\to{\rm sech}(z)$ and $K(m)\to\infty$ as $m\to 1^-$ 
and hence \eqref{cnoidal1} reduces to the well-known formula for the KdV soliton.}.

\begin{figure}[htbp]
\begin{center}
(a)\includegraphics[scale=0.75]{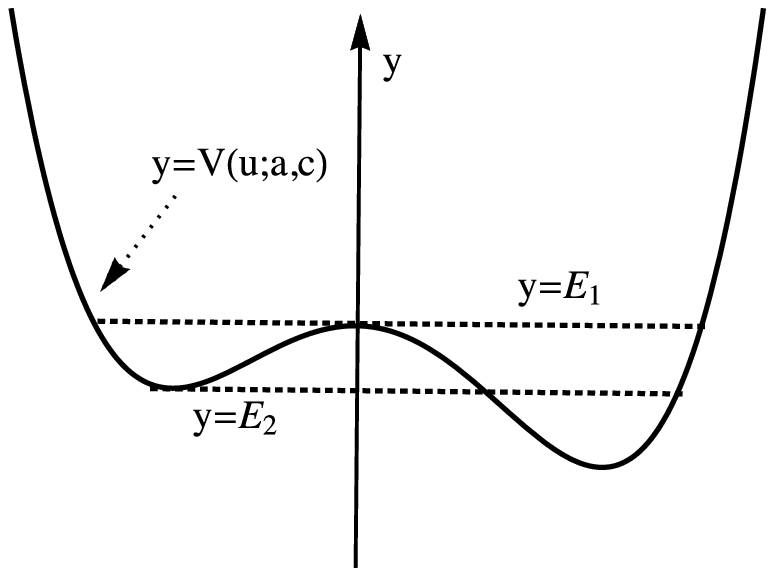}\quad(b)\includegraphics[scale=0.75]{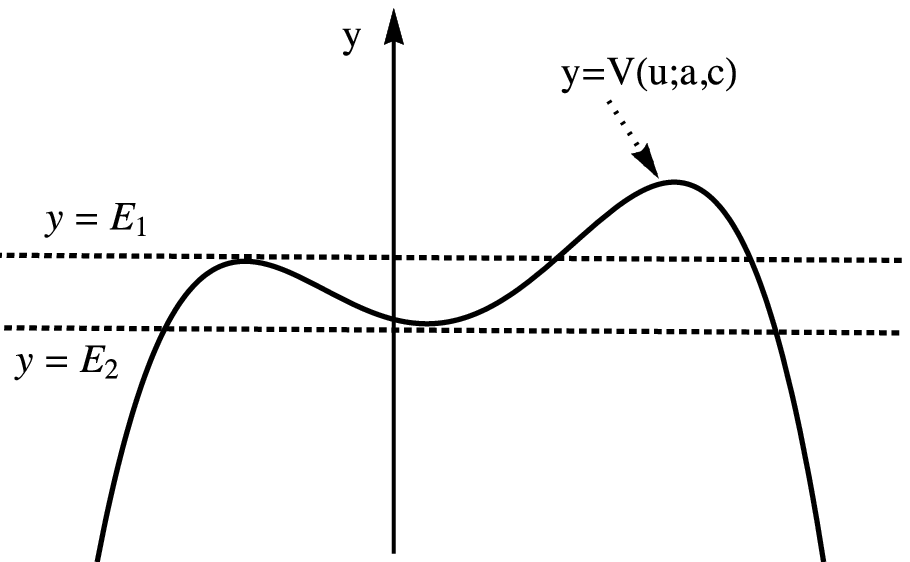}
\end{center}
\caption{Representative graphics of the effective potential energy $V(\cdot\,;a,c)$ 
in the case of the (a) focusing and (b) defocusing modified KdV equations. 
In the focusing case, periodic solutions exist only for $E>E_2(a,c)$. 
There are two families of periodic solutions when $E\in(E_2,E_1)$, 
while there is only one family when $E>E_1$.  
In the defocusing case, periodic solutions exist only for $E\in(E_2,E_1)$
and there is only one family of solutions.
}\label{fig1}
\end{figure}

Another example of explicit periodic traveling waves is when $f(u)=\pm \frac{1}{3} u^3$ in \eqref{gkdv},
which are called the focusing (+) and defocusing (-) modified KdV equations. 
In the focusing case, $V(\cdot\,;a,c)$ describes a bistable double-well potential 
and hence simple phase plane analysis implies that 
non-constant periodic solutions of \eqref{quad2} exist, depending on $a$, $E$, and $c$, 
provided that the quartic polynomial $E-V(\cdot\,;a,c)$ has four distinct roots 
(i.e. its discriminant is non-zero), say, $\alpha$, $\beta$, $\gamma$ and $\delta$, 
with at least two of them being real; see Figure~\ref{fig1}, for instance. 
The general solution is rather tedious with a number of cases to consider, 
but it amounts to reducing the elliptic integral to standard form. In fact, if we write the solution in the form
\[
\int \frac{du}{\sqrt{(u-\alpha)(u-\beta)(u-\gamma)(u-\delta)}}= \omega t
\]
for some appropriate constant $\omega\in\mathbb{R}$ and make the M\"obius tranformation 
\[
v = \frac{au+b}{u-\alpha}
\]
for $a, b$ constants, 
which maps one of the roots of the quartic to $\infty$, then the above elliptic integral  becomes 
\[
\int \frac{dv}{\sqrt{P_3(v)}}= \omega t   
\]
for some cubic polynomial $P_3$, which may be analyzed as above. 
Clearly, the root $\alpha$ must not be an end point of the interval over which the solution oscillates.

In the case of $a=0$, in particular, the effective potential $V(\cdot\,;0,c)$ is symmetric (even) and
the expressions for explicit solutions simplify greatly.  
Specifically, there are two families of solutions, one in terms of the Jacobi elliptic cnoidal function ${\rm cn}$ 
and the other in terms of the Jacobi dnoidal function ${\rm dn}$.
The cnoidal waves correspond to solutions with $E>0$ and their orbits in phase space 
lie outside the separatrix; see Figure \ref{fig1}(a). 
The dnoidal waves correspond to solutions with $E<0$
and their orbits in phase space are bounded by a homoclinic orbit (corresponding to a solitary wave). 
In the case of $c<0$ and $E<0$, for instance, the quartic on the right side of the profile ODE
\[
u_z^2 = 2 E - c u^2 - \frac{u^4}{6}
\]
has four real roots, $u=\pm k_1,u=\pm k_2$, 
with $k_1=-3(c+\sqrt{c^2+ \frac{4E}{3}})$ and $k_2=-3(c-\sqrt{c^2+ \frac{4E}{3}})$. 
The change of variables $u = k_2 v$ gives 
\[
v_z^2 = \frac{k_2^2}{6}(v^2-1)\Big(\frac{k_2^2}{k_1^2}-v^2\Big).
\]
Note that the elliptic function $y(z)=\mbox{dn}(z; k)$ solves 
\[
y_z^2 = (y^2-1)(1-k^2-y^2).
\]
Therefore $v = \mbox{dn}(\frac{k_2}{\sqrt{6}}z;k)$, $k^2 = 1 - \frac{k_1^2}{k_2^2}$. 
The other case, when $c<0$ and $E>0$, 
which is solvable in terms of the cnoidal function, follows similarly.

In the defocusing case, non-constant periodic solutions exist, depending on $a$, $E$, and $c$, 
only when the quartic polynomial $E-V(\cdot\,;a,c)$ has four distinct real roots; see Figure~\ref{fig1}(b). 
In the case $a=0$, the potential is likewise symmetric and one may find explicit solutions
in terms of the Jacobi snoidal function.

\subsection{General Existence Theory}

While an extensive amount of literature is devoted to studying 
properties of explicit, periodic traveling waves of \eqref{gkdv}, 
e.g. the elliptic function solutions of the KdV and modified KdV equations discussed above, 
our modulational instability theory does not make use of such formulae.
(Nevertheless, we will strive to state our general results in the context of explicit solutions when possible.)
Rather, of importance to us are properties of the \emph{parameterization}
of periodic traveling waves with respect to various coordinate systems and the Jacobians between them.
As a first step in the direction, we remark that 
the Jacobian of the map \eqref{trans} from $\mathbb{R}^3$ to itself is non-singular, 
provided that the cubic polynomial $E-V(\cdot\,;a,c)$
has three distinct roots, a property fundamentally related to the sign of its \emph{discriminant}.  
As a result, in the set where the discriminant of the polynomial is positive, $E-V(\cdot\,;a,c)$ has three distinct real roots
and the transformation \eqref{trans} be smoothly invertible.
Consequently, the cnoidal wave solutions \eqref{cnoidal1} of the KdV can be considered 
as a four-parameter family of periodic traveling waves
parametrized by the ``natural" ODE parameters $a$, $E$, $c$, and $z_0$, 
whose profile and fundamental period depends smoothly
on these parameters.  Analogous results hold for elliptic function solutions 
of the focusing and defocusing modified KdV equations: within the domain
of existence they can be smoothly parametrized in terms of $E$, $c$, and $z_0$ 
(recalling that $a=0$ in these explicit solutions).

With the above in mind, we note that if $f$ in \eqref{gkdv} is a polynomial nonlinearity 
then the zero set of the discriminant of $E-V(\cdot;a,c)$, denoted as
\[
\Gamma:=\Big\{(a,E,c)\in\mathbb{R}^3:\textrm{disc}\Big(E+au-\frac{1}{2}cu^2-F(u)\Big)=0\Big\},
\]
defines a \emph{variety} dividing the parameter space into open sets with a constant number of periodic solutions.  Solutions of \eqref{quad2}
with $(a,E,c)\in\Gamma$ correspond to either equilibrium or homoconic solutions of \eqref{quad2}, 
which, respectively, correspond to constant solutions and solitary waves of \eqref{gkdv}.
In particular, $(0,0,c)\in\Gamma$ represents the solitary waves homoclinic to zero.  

\begin{figure}
\begin{center}
\includegraphics[scale=0.85]{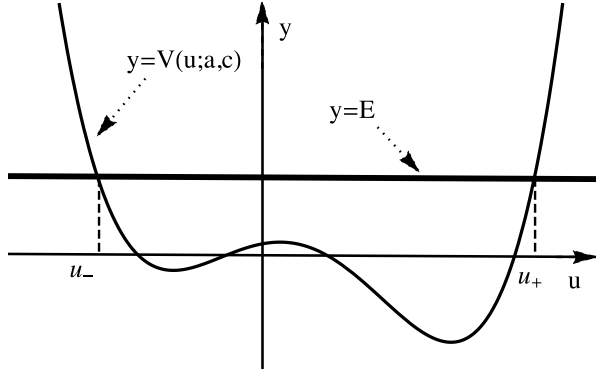}
\end{center}
\caption{Schematic drawing showing the locations of the roots $u_+$ and $u_-$.
}\label{fig2}
\end{figure}

In order to ensure the existence of periodic orbits of \eqref{quad2}, we take $(a,E,c)\in\mathbb{R}^3\setminus\Gamma$
so that there exist simple and distinct roots of the equation $E=V(\cdot;a,c)$ and that there are real roots $u_{\pm}$ satisfying $u_-<u_+$
such that $V(u;a,c)<E$ for $u\in(u_-,u_+)$; see Figure \ref{fig2}.  It follows that the roots $u_{\pm}$ are smooth functions of the traveling wave parameters $(a,E,c)$ and,
without loss of generality, we can assume $u(0)=u_-$.  It follows that the corresponding periodic solution of \eqref{quad2} is an  even function of $z=x-ct$
and that its period can be expressed via quadrature as
\[
T=T(a,E,c)=\sqrt{2}\int_{u_-}^{u_+}\frac{du}{\sqrt{E-V(u;a,c)}}=\frac{\sqrt{2}}{2}\oint_\gamma\frac{du}{\sqrt{E-V(u;a,c)}},
\]
where integration over $\gamma$ represents a complete integration from $u_-$ to $u_+$ and then back to $u_-$ again. Notice however
that the branch of the square root must be chosen appropriately in each direction.  Alternatively, the contour $\gamma$
can be interpreted as a closed loop (Jordan curve) in the complex plane that encloses a bounded set containing both $u_{\pm}$.  By a standard
procedure, the above integral can be regularized at the square root branch points and hence represents a $C^1$ function of $a$, $E$, and $c$.

In general, \eqref{gkdv} admits three conserved quantities, 
physically corresponding to the Hamiltonian, the mass, and the momentum.
Given a $T$-periodic traveling wave $u(\cdot;a,E,c)$ of \eqref{gkdv}, they are defined by
\begin{equation}\label{cons1}
\begin{aligned}
H(a,E,c)&=\int_0^T\Big(\frac{u_x^2}{2}-F(u)\Big)~dx=\frac{\sqrt{2}}{2}\oint_\gamma\frac{E-V(u;a,c)-F(u)}{\sqrt{E-V(u;a,c)}}~du,\\
M(a,E,c)&=\int_0^Tu~dx=\frac{\sqrt{2}}{2}\oint_\gamma\frac{u~du}{\sqrt{E-V(u;a,c)}},\\
P(a,E,c)&=\int_0^Tu^2~dx=\frac{\sqrt{2}}{2}\oint_\gamma\frac{u^2~du}{\sqrt{E-V(u;a,c)}},
\end{aligned}
\end{equation}
respectively, where the integral over $\gamma$ is defined as before.  As above, these integrals can be regularized at their square root branch points and hence
represent $C^1$ functions of the traveling wave parameters.  As we will see, the gradients of the period and the above conserved quantities plan a
fundamental role in the stability of the periodic traveling wave solutions of \eqref{gkdv}.

In summary, \eqref{gkdv} in the case of polynomial nonlinearity admits 
a four-parameter family of periodic traveling waves, 
smoothly parametrized by $a$, $E$, $c$ in \eqref{quad2} and an auxiliary parameter $z_0$
associated with the translational invariance of the equation.
Solitary waves, whose profile asymptotically vanishes, correspond to $a=E=0$, 
and hence they form a co-dimension two subset of the full set of traveling waves.

\section{Formal Asymtotics and Whitham's Theory of Modulations}\label{s:whitham}

We are interested in the stability of periodic traveling waves of nonlinear dispersive equations
to slow modulations. That means, we wish to understand 
the long-term dynamics of ``slightly" modulated periodic traveling waves, in particular, to determine 
whether the long-term dynamics will be appropriately described by small modulations of the carrier wave.
In this section, we illustrate some asymptotic calculations 
that shed light on the long-term dynamics and hence on the modulational stability.
Mainly, we will utilize a formal multiple-scales expansion
to derive a first order system of equations governing 
the slow dynamics of averaged quantities of the modulated wave over the fast scale.
Such an approach was pioneered by Whitham in the 1960's and early 1970's 
(see \cite{Whitham.1965}, \cite{Whitham.1967}, \cite{Whitham.book}, for instance), 
and it has since been widely applied in both the physics and mathematics communities. 
As such, the formal theory has been affectionately termed Whitham's theory of modulations, 
or simply, Whitham's theory.

To illustrate the key ideas, we begin by considering slow modulations of a linear dispersive wave packet, where the modulational dynamics are simplified due to the direct coupling of spatial and temporal frequencies via the linear dispersion relation.  We will then extend this calculation to nonlinear wave packets where one must make use of a nonlinear dispersion relation to close the modulation system.  Section \ref{s:validate}, we will discuss rigorously validating the predictions from these calculations.

\subsection{Linear Dispersive Waves}
Notice that any constant function $u(x,t)=u_0$ is a solution of \eqref{gkdv}.
Seeking nearby solutions of the form 
\[
u(x,t)=u_0+\varepsilon v(x,t), \qquad |\varepsilon|\ll 1,\quad v=\mathcal{O}(1),
\]
we arrive at that 
\[
v_t+v_{xxx}+f'(u_0)v_x=\mathcal{O}(\varepsilon).
\]
Ignoring all $\mathcal{O}(\varepsilon)$, we then find that as long as $v$ remains $\mathcal{O}(1)$ 
we can approximate $v$ by a solution of the linear equation
\begin{equation}\label{lingKdV}
v_t+v_{xxx}+f'(u_0)v_x=0.
\end{equation}
It admits plane wave solutions of the form
\begin{equation}\label{linwave}
v(x,t)=Ae^{i(kx-\omega t)}+c.c.,
\end{equation}
where $A$ is the (real) amplitude, $k$ is the wave number, 
and the temporal frequency $\omega$ is related to $u_0$ and $k$ via the dispersion relation
\begin{equation}\label{disp}
\omega=f'(u_0)k-k^3.
\end{equation}

We now wish to describe the modulation of the linear wave packets. 
That is, we consider solutions of \eqref{lingKdV}
of the form \eqref{linwave}, where the parameters $A$, $k$, and $\omega$ are allowed to depend on space and time.
First, note that if $A$, $k$, and $\omega$ are allowed to vary significantly over one wavelength of the ``carrier" wave \eqref{linwave}, taken
with $A$, $k$, and $\omega$ constant, then there is no reason to expect that the evolving wave packet should be described
as a modulated periodic wave.  Hence, we only consider modulations, where $A$, $k$, and $\omega$ vary slowly
on space and time.  This naturally leads one to consider two separate space and time
scales, namely the ``fast" oscillations of the carrier wave and the slow evolution of $A$, $k$, and $\omega$.  

Utilizing a multiple scales approach, we thus seek a modulated plane wave solution of \eqref{lingKdV} of the form
\begin{equation}\label{lin_mod}
v(x,t,X,S)=A(X,S)e^{i(k(X,S)x-\omega(X,S)t)}+c.c.,
\end{equation}
where $(x,t)$ denote the ``fast" variables and $(X,S)=(\varepsilon x,\varepsilon t)$ denote the ``slow" 
variables. 
Here and elsewhere, $\varepsilon$ refers to a small, but non-zero, quantity.  
Substituting this into \eqref{lingKdV} and collecting powers of $\varepsilon$ we find that
\begin{multline}\label{exp1}
(-i\omega+ikf'(u_0)-ik^3)+i\varepsilon A\Big(\Big(\frac{\partial k}{\partial S}+(u_0-3k^2)\frac{\partial k}{\partial X}\Big)(x-(u_0-3k^2 t))\Big)\\
+\varepsilon\Big(\frac{\partial A}{\partial S}+(u_0-3k^2)\frac{\partial A}{\partial X}-3kA\frac{\partial k}{\partial X}\Big)=0.
\end{multline}
The $\mathcal{O}(1)$ equation simply requires that the slowly evolving temporal and spatial frequencies satisfy the dispersion relation \eqref{disp}, i.e.
\[
\omega(k)=f'(u_0)k-k^3,\qquad k=k(X,S),
\]
for all $X,S\in\mathbb{R}$.  Since this forces $\omega$ to be real valued regardless of the slow evolution of the wave number $k$, we see
that the modulated solutions \eqref{lin_mod}  are bounded for all time, and hence the constant solution $u_0$ is considered to be 
modulationally stable.

To describe the slow evolution more precisely, note that at 
$\mathcal{O}(\varepsilon)$ the expansion \eqref{exp1} yields two equations for the slow dynamics of $k$ and $A$:
\begin{equation}\label{w1}
\begin{aligned}
\frac{\partial k}{\partial S}+\omega'(k)\frac{\partial k}{\partial X}=0,\qquad
\frac{\partial A}{\partial S}+(\omega'(k)A)_X=0, 
\end{aligned}
\end{equation}
where $\omega'(k)=f'(u_0)-3k^2$ depends on $X$ and $S$ through the evolution of $k$.  
Together with the dispersion relation, the system \eqref{w1} describes, to leading order, the evolution of slowly modulated wave packets in \eqref{lingKdV}.
In particular, since $\omega(k)$ is real-valued for all $k\in\mathbb{R}$, the linearization of this system about $(k,A)=(k_0,A_0)$, corresponding a linear wave train
of the form \eqref{linwave} with $k=k_0$ and $A=A_0$ constant, is hyperbolic, and hence both $k$ and $A$ will
remain bounded for all time.  Again, this implies the modulational stability of the constant solution.

Note that the Whitham modulation equations \eqref{w1} for slow modulations of a constant state prescribe the evolution of only \emph{two}
quantities, namely the wave number and the amplitude of the wave. Moreover, the evolution equation
for $k$ is \emph{independent of the evolution of the amplitude}.  These observations are due to the fact that
the dispersion relation \ref{disp} directly relates the temporal and spatial frequencies.  

\subsection{Nonlinear Dispersive Waves}
We now wish to understand the long-time dynamics of an arbitrary amplitude, slowly modulated  periodic traveling waves of \eqref{gkdv}.
The main technical difference is that the spatial and temporal frequencies are no longer coupled through
the linaer dispersion relation, and hence we expect to find a modulation system describing the slow evolution of \emph{three} of the characteristic parameters of the wave.

Recall from Section \ref{s:exist} that periodic traveling waves of \eqref{gkdv} 
are parametrized (up to translations) by three parameters $a$, $E$, and $c$, 
and the period depends on the parameters.
While this parametrization is natural from the existence standpoint, we will see below that it is
is not well-suited to the formal asymptotic description of Whitham's theory of modulations.
(We utilize this parametrization in the rigorous theory in the following sections, however.) 
Rather, the modulation equations are most naturally described
in terms of the slow evolution of the period and particular conserved quantities associated with the modulated solution.  The Jacobian
between these parametrizations will play a vital role in our theory.

We seek an asymptotic description of modulated periodic traveling waves of \eqref{gkdv}.
To account for the separate fast and slow scales of space and time involved in such modulations,  
we utilize a multiple scales expansion.
We introduce ``slow" variables $(X,S)=(\varepsilon x, \varepsilon t)$ 
and note that, in the slow coordinates, \eqref{gkdv} becomes
\begin{equation}\label{slowgkdv}
u_S=\varepsilon u_{XXX}+f(u)_X.
\end{equation}
Following Whitham \cite{Whitham.1965,Whitham.book}, we seek a solution of \eqref{slowgkdv} of the form
\begin{equation}\label{wkb}
u(X,S)
=u^{(0)}\Big(\frac{\phi(X,S)}{\varepsilon},X,S\Big)
+\varepsilon u^{(1)}\Big(\frac{\phi(X,S)}{\varepsilon},X,S\Big)+\mathcal{O}(\varepsilon^2),
\end{equation}
where the phase $\phi$ is chosen to ensure that $u^{(j)}$, $j=0,1$, are 1-periodic functions of the variable $y=\phi(X,S)/\varepsilon$. In particular,
we note that the local period of oscillation of the modulated wave train $u^0$ is $\varepsilon/\partial_X\phi$, where we assume the unknown phase
\emph{a priori} satisfies that $\partial_X\phi\neq 0$.  Notice that periodic traveling waves of \eqref{gkdv} fit this asymptotic description
by taking $u=u^{(0)}$ constant in the slow variables and choosing $\phi(X,S)=k(X-cS)$.  That is,
periodic traveling waves of \eqref{gkdv} correspond to \emph{constant solutions}  in the slow coordinates.

We now substitute \eqref{wkb} into \eqref{slowgkdv} and collect like powers
of $\varepsilon$, resulting in a hierarchy of equations that must be satisfied.  At the order of $\mathcal{O}(\varepsilon^{-1})$, we find that
\[
\phi_S\partial_yu^{(0)}=(\phi_X\partial_y)^3u^{(0)}+(\phi_X\partial_y)f(u^{(0)}) 
\]
which is simply \eqref{quad2} for $u^{(0)}$ differentiated in the variable $y$,  
under the identification $k=\phi_X$ as the spatial frequency and $c=\omega/k$ as the wave speed, where $\omega=-\phi_S$
is the temporal frequency of the modulation, all of which are considered functions of the slow variables $X$ and $S$.  For fixed $X$ and $S$, thus we may
choose $u^{(0)}$ to be a periodic traveling wave of \eqref{gkdv}, and hence of the form
\begin{equation}\label{profile0}
u^{(0)}(y,X,S)=\bar{u}(y;a(X,S),E(X,S),c(X,S))
\end{equation}
for some solution $\bar{u}$ of \eqref{quad2}, where $a$, $E$, and $c$ are independent of $y$.
In particular, the compatibility condition $\phi_{SX}=\phi_{XS}$ implies that
\begin{equation}\label{w2}
k_S+\omega_X=0,
\end{equation}
which serves as a nonlinear dispersion relation. 
It effectively replaces the dispersion relation for linear waves.  
The linear case satisfy \eqref{w2} thanks to the dispersion relation and the first equation in \eqref{w1}.
In the literature, \eqref{w2} is sometimes referred to as the equation
of conservation of waves.  

At the order of $\mathcal{O}(1)$, continuing, 
\begin{equation}\label{oeps}
\partial_z\mathcal{L}[u^{(0)}]u^{(1)}=\partial_Su^{(0)}-\partial_Xf(u^{(0)})-\partial_X(\phi_X^2\partial_y^2)u^{(0)},
\end{equation}
supplemented with $2\pi$-periodic boundary conditions, where 
\[
\mathcal{L}[u^{(0)}]=\partial_z^2-su^{(0)}+f'(u^{(0)}), \qquad z=k y,
\]
is recognized as the linear operator obtained from linearizing \eqref{gkdv}
about the carrier wave $u^{(0)}$.   The Fredholm alternative implies that \eqref{oeps} has a solution, provided that
the right hand side is orthogonal to the kernel of the adjoint operator $(\partial_y\mathcal{L}[u^{(0)}])^\dag=-\mathcal{L}[u^{(0)}]\partial_y$ 
in the sense of $L^2(0,2\pi)$.  Since $u^{(0)}$ is not
constant, it can readily be shown by differentiating \eqref{quad2} that 
\[
{\rm ker}(\mathcal{L}[u^{(0)}]\partial_y)={\rm span}\{1,u^{(0)}\}.
\]
Hence \eqref{oeps} has a solution, provided that the orthogonality conditions
\begin{equation}\label{w3}
\left\{\begin{aligned}
&M_S+\partial_X G(u^{(0)})=0,\\
&P_S+\partial_XQ(u^{(0)})=0
\end{aligned}\right.
\end{equation}
hold, where 
\[
M(X,S)=\int_0^{1}u_0(y,X,S)dy\quad\text{and}\quad P(X,S)=\int_0^{1}u_0^2(y,X,S)dy
\]
are functions of the slow variables $X$ and $S$, 
and the fluxes $G$ and $Q$ denote the inner product of the right hand side of \eqref{oeps} with $1$ and $u_0$, respectively;
in particular, they are $y$-independent functions of $X$ and $S$ that may be written out explicitly, if desired;  see \cite{JZ}, for instance.

Together, three equations in \eqref{w2} and \eqref{w3} form the Whitham averaged system for \eqref{gkdv}, describing the mean behavior
of the slowly varying functions $a$, $E$, and $c$ in \eqref{profile0} \emph{implicitly} in terms of the slow evolution
of the functions $k$,  $M$, and $P$.  In particular, it formally follows that 
\eqref{gkdv} has a modulated wave train solution that to leading order can
be described as
\[
u(x)=\bar{u}\Big(\frac{\phi(X,S)}{\varepsilon};a(X,S),E(X,S),c(X,S)\Big)+\mathcal{O}(\varepsilon),
\]
where $a$, $E$, and $c$ evolve in the slow variables in such a way that 
$k$, $M$, and $P$ evolve according to the (nonlinear) first order system \eqref{w2}-\eqref{w3}.

To ascertain the modulational stability of periodic traveling wave $u_0(a_0,E_0,c_0)$ of \eqref{gkdv},
we seek to understand whether the functions $a$, $E$, and $c$ in \eqref{profile0} evolve away from $(a_0,E_0,c_0)$
in the slow scale or not.  The modulation system \eqref{w2}-\eqref{w3}, however, describes the evolution of $a$, $E$ and $c$
\emph{implicitly} through the evolution of $k$, $M$, and $P$.  As such, if we want to measure deviation
of $(a,E,c)$ from $(a_0,E_0,c_0)$ in terms of the evolution of $(k,M,P)$, we must require that, at least in the vicinity of the carrier wave $u_0$,
nearby periodic traveling waves of \eqref{gkdv} can be \emph{re-parametrized} by the parameters $k$, $M$, and $P$.
By the implicit function theorem, such a (local) re-parametrization is possible, provided that the Jacobian matrix
\begin{equation}\label{jacob1}
\frac{\partial(k,M,P)}{\partial(a,E,c)}
\end{equation}
is non-singular at $(a_0,E_0,c_0)$.

Under the assumption that the matrix \eqref{jacob1} is non-singular, the modulational stability 
can be analyzed by linearizing the modulation equations
\eqref{w2}-\eqref{w3} about the constant solution $(a_0,E_0,c_0)$.  We seek solutions
of \eqref{w2}-\eqref{w3} of the form
\[
(a(X,S),E(X,S),c(X,S))=(a_0,E_0,c_0)+\varepsilon(\tilde a(X,S),\tilde E(X,S),\tilde c(X,S)),
\]
which, upon substituting into \eqref{w2}-\eqref{w3} and dropping all $\mathcal{O}(\varepsilon^2)$
terms\footnote{Note the $\mathcal{O}(1)$ terms cancel out due to the definition of $(a_0,E_0,c_0)$, so the leading order term is $\mathcal{O}(\varepsilon)$.}, yields the linear system
\begin{equation}\label{linwhitham}
\Big(\frac{\partial( k_0,M_0, P_0)}{\partial(a,E,c)}\partial_S(\tilde a,\tilde E,\tilde c)+\frac{\partial( \omega_0, G_0, Q_0)}{\partial(a,E,c)}\partial_X(\tilde a,\tilde E,\tilde c)\Big)=0,
\end{equation}
where the $0$ subscript denotes evaluation at $(a_0,E_0,c_0)$.
This is a constant coefficient, linear system of PDE's, and we may seek solutions
of the form 
\[
(\tilde a,\tilde E,\tilde c)(X,S)=e^{\lambda S-ik X}(\tilde a_0,\tilde E_0,\tilde c_0),
\qquad\lambda\in\mathbb{C},\quad k\in\mathbb{R}
\]
for some $(\tilde a_0,\tilde E_0,\tilde c_0)\in\mathbb{R}^3$,
in which case the temporal and frequencies of the perturbation can be determined from the dispersion relation
\[
\mathcal{D}_W(\lambda,k):=\det\left(\lambda\frac{\partial(k_0,M_0,P_0)}{\partial(a,E,c)}-ik\frac{\partial(\omega_0,G_0,Q_0)}{\partial(a,E,c)}\right)=0,
\]
a homogeneous cubic polynomial in the variables $\mu$ and $k$.  In particular, $\mathcal{D}_W(\mu,k)=\mu^3\mathcal{D}_W(1,k/\mu)$
and hence if $\{z_j\}_{j=1}^3$ are the roots of  
\begin{equation}\label{Wdisp}
\widehat{\mathcal{D}}_W(z):=\mathcal{D}_W(1,z)
\end{equation}
then the dispersion relation defines three spectral
curves
\[
\lambda_j(k)=iz_jk,\qquad j=1,2,3.
\]
In particular, $u_0(a_0,E_0,c_0)$ should be modulationally stable, provided that the reduced dispersion relation 
has three real roots, while $u_0$ should be modualtionally unstable if it admits at least one complex root with non-zero imaginary part.

The above formal analysis suggests that a necessary condition for the modulational stability of a periodic wave train of \eqref{gkdv} is
for the linearized Whitham system \eqref{linwhitham} to be hyperbolic, i.e. the matrix ${\displaystyle \frac{\partial(k_0,M_0,P_0)}{\partial(\omega_0,G_0,Q_0)}}$
has all real eigenvalues.  In the next section, we rigorously derive a modulational stability theory whose main result, Theorem \ref{t:whitham},
relates the stability of a periodic wave train of \eqref{gkdv} to slow modulations to the eigenvalues of a $3\times 3$ matrix.  These eigenvalues
end up agreeing precisely with the roots of $\widehat{\mathcal{D}}_W$, providing a rigorous justification of the above formal calculations; see Theorem \ref{t:validate} below.

\section{Rigorous Theory of Modulational Instability}\label{s:theory}

We discuss how to make rigorous the formal arguments sketched in the previous section in various settings.
In the process, we provide a general methodology 
of detecting modulational instabilities for a large class of nonlinear dispersive equations.

\subsection{Analytic Setup}

Let $u_0(\cdot;a,E,c)$ denotes a $T=T(a,E,c)$ periodic traveling wave solution of the generalized gKdV equation \eqref{gkdv}.  The idea 
of ``stability" can formally be stated as requiring that if a solution to \eqref{gkdv} starts close to the solution $u_0$ then it stays close to the solution
$u_0$ for all time.  To understand the dynamics of solutions that start near $u_0$, one typically begins by linearizing the governing
PDE about the solution $u_0$, a process we now describe.

First, notice that $u_0$ is a stationary $T$-periodic solution of the PDE
\begin{equation}\label{travelgkdv}
u_t=cu_z+u_{zzz}+f(u)_z,
\end{equation}
which is simply the gKdV equation \eqref{gkdv} written in the traveling spatial frame $z=x-ct$.  We now seek a solution
$u$ of \eqref{travelgkdv} with initial data $u(0)=u_0+v_0$ for some  $v_0\in H^3(\mathbb{R})$ with $\|v_0\|_{H^3(\mathbb{R})}$
sufficiently small\footnote{The requirement that the perturbation be integrable is not the only choice possible.  Other natural
candidates of classes of perturbations include periodic perturbations with fundamental period $nT$ for some $n=1,2,3,\ldots$.  While
such periodic classes of perturbations are natural from a variationsl viewpoint, they may impose artificial constraints on the physical
problem.  As we will see below, the class of localized (i.e. integrable on the line) perturbations includes information about all quasi-periodic perturbations.}. 
Writing $v(t)=u_0-u(t)$ for $t>0$ for which $u(t)$ is defined\footnote{Here, we ignore issues
of well-posedness of \eqref{gkdv} to such initial data.}, 
it follows that $v$ satisfies the evolution equation
\begin{equation}\label{evol2}
v_t=Lv+\mathcal{N}[u_0,v],
\end{equation}
where here $L=\partial_u F(u_0):H^3(\mathbb{R})\subset L^2(\mathbb{R})\to L^2(\mathbb{R})$ is the linear operator
\[
L=\partial_z\left(\partial_z^2+c+f'(u_0)\right),
\]
and where $\mathcal{N}[u_0,v]=\partial_z\left(f(u_0+v)-f(u_0)-f'(u_0)v\right)$ is a nonlinear operator satisfying the quadratic estimate
\[
\left\|\mathcal{N}[u_0,v]\right\|_{L^2(\mathbb{R})}\leq C\|v\|_{H^1(\mathbb{R})}^2
\]
for some constant $C=C(\|u_0\|_{H^1(\mathbb{R})})>0$ independent of $v$.  
We expect that as long as $v_0$ is sufficiently small in, say, $L^2(\mathbb{R})$ 
then the solution $v(t)$ should be small at least for small, again in $L^2(\mathbb{R})$ say, at least for short times.  Moreover, so long as $v(t)$ remains
sufficiently small it is reasonable to expect that the nonlinear evolution \eqref{evol2} may be well approximated by the linear evolution equation
\begin{equation}\label{linevol}
v_t=Lv
\end{equation}
posed on $H^3(\mathbb{R})\subset L^2(\mathbb{R})$.

Due to linearity, we may of course relax the smallness assumption in the analysis of \eqref{linevol}.  
The underlying wave $u_0$ is said to be \emph{linearly stable} (to perturbations in $L^2(\mathbb{R})$) provided that all
solutions of \eqref{linevol} remain bounded in time, while it is \emph{linearly unstable} (to $L^2(\mathbb{R})$ perturbations)
if there exists some initial data $v_0\in H^3(\mathbb{R})$ such that the associated solution of \eqref{linevol}
grows in time.

The linear stability of $u_0$ may be further reduced to a problem in spectral analysis as follows.
As the linear evolution \eqref{linevol} is autonomous in time, we may take the Laplace transform in time, resulting in the problem
\begin{equation}\label{spec}
Lv=\lambda v,
\end{equation}
where here $v(z)=v(z;\lambda)$ denotes the Laplace transform of the function $v(t,z)$ in \eqref{linevol} and $\lambda\in\mathbb{C}$ is the Laplace (temporal) frequency.  
Equation \eqref{spec} is a spectral problem posed on $L^2(\mathbb{R})$, and as such the stability of $u_0$ may be studied by analyzing the spectrum
of $L$ considered as a densely defined operator on $L^2(\mathbb{R})$.
Indeed, if the $L^2(\mathbb{R})$-spectrum of $L$ intersects the open right half plane in $\mathbb{C}$,
then one expects the solution of \eqref{spec} to grow in time and hence the underlying solution $u_0$ should be unstable.  With this in mind,
we say that the underlying wave $u_0$ is \emph{spectrally unstable} (to $L^2(\mathbb{R})$ perturbations) if there exists a $\lambda$ with $\Re(\lambda)>0$ in the
 $L^2(\mathbb{R})$- spectrum of $L$.  Otherwise, we say $u_0$ is \emph{spectrally stable} (to $L^2(\mathbb{R})$ perturbations)\footnote{Generally, spectral stability does \emph{not}
imply linear stability, as is familiar from basic ODE theory.  Nevertheless, spectral \emph{instability} often does imply linear (and nonlinear) instability.}.
Notice that since the spectral problem \eqref{spec} is left invariant under the transformations
\[
v\mapsto \bar{v}\quad\textrm{and}\quad\lambda\mapsto\bar{\lambda},
\]
where here the ``bar" denotes complex conjugation, 
and
\[
z\mapsto -z\quad\textrm{and}\quad\lambda\mapsto -\lambda
\]
the spectrum of $L$ is symmetric with respect to reflections about both the real and imaginary axes.  Consequently, $u_0$ is spectrally stable
if and only if the $L^2(\mathbb{R})$-spectrum of $L$ is confined to the imaginary axis.

As \eqref{spec} is a spectral problem for a linear operator with periodic coefficients, its spectrum is most conveniently characterized
via Floquet-Bloch theory.  Standard results from this theory \cite{Chicone.book,RS4} dictate that 
non-trivial solutions of \eqref{spec} can not be integrable on $\mathbb{R}$, i.e. they can not have finite norm in $L^p(\mathbb{R})$ for any $1\leq p<\infty$, 
and that they can at best be bounded functions on the line.  Further, any bounded solution of \eqref{spec} must be of the form
\[
v(x)=e^{i\xi x}w(x)
\]
for some $w\in L^2_{\rm per}([0,T])$ and $\xi\in[-\pi/T,\pi/T)$.  In particular, a given $\lambda\in\mathbb{C}$
belongs to the $L^2(\mathbb{R})$-spectrum of the operator $L$ if and only if there is a non-trivial solution
of the quasi-periodic spectral problem
\[
\left\{\begin{aligned}
&Lv=\lambda v,\\
&v(x+T)=e^{i\xi T}v(x)
\end{aligned}\right.
\]
for some $\xi\in[-\pi/T,\pi/T)$.  Alternatively, $\lambda\in\mathbb{C}$ if and only if there exists a $\xi\in [-\pi/T,\pi/T]$ such
that there is a non-trivial $T$-periodic solution of the equation
\[
L_\xi w=\lambda w,\quad \textrm{where}\quad L_\xi:=e^{-i\xi x}Le^{i\xi x}.
\]
For each $\xi\in[-\pi/T,\pi/T)$, the spectrum of the operator $L_\xi$, considered as a densely defined operator on $L^2_{\rm per}([0,T])$, is comprised merely of discrete eigenvalues with finite multiplicity and, furthermore,
\[
{\rm spec}_{L^2(\mathbb{R})}(L)=\bigcup_{\xi\in[-\pi/T,\pi/T)}{\rm spec}_{L^2_{\rm per}([0,T])}(L_\xi).
\]
The parameter $\xi$ is referred to as the Bloch frequency and the operators $L_\xi$ are called the Bloch operators associated to $L$.

Rather than analyzing the essential spectrum of the operator $L$ directly, the above discussion indicates that 
we may rather choose to study the $T$-periodic eigenvalues of the Bloch operators
$L_\xi$ for each $\xi\in[-\pi/T,\pi/T)$.  From a practical level, one does not expect to be able to explicitly compute the eigenvalues
of $L_\xi$ for an arbitrary Bloch frequency $\xi$, except outside of a few very special cases\footnote{For example, completely integrable
systems \cite{Bottman.Deconinck.2009}.}.  Thankfully, however,for the purposes of modulational stability analysis
we need only consider the spectrum of the operators $L_\xi$ in a neighborhood
of the origin in the spectral plane, and only for $|\xi|\ll 1$.  
Indeed, notice that if the eigenvalues of $L_0$ are confined to the imaginary axis,
then $u_0$ is spectrally stable to perturbations which are themselves $T$-periodic.  
When $0<|\xi|\ll 1$, the period of the perturbation is neary  that of the 
underlying wave $u_0$, and hence the effect of the perturbation is seen only on very large space and time scales.  Hence, the spectrum
of the operators $L_\xi$ with $0<|\xi|\ll 1$ describes the stability of $u_0$ to \emph{long wavelength perturbations}.  

Modulational perturbations clearly fall into the catagory of long-wavelength perturbations, but they form a special subclass where the effect of the perturbation is simply
to modulate the defining characteristics of the wave, here given by the parameters $a$, $E$, $c$, and the translational mode $z_0$.  The key point
is that variations in these defining parameters provide spectral information of $L_0$ at $\lambda=0$ in the spectral plane.  Formally, this can
be seen by differentiating the profile equation \eqref{quad1} with respect to $a$, $E$, $c$, and $z$ and noting that, 
\begin{equation}\label{kerform}
L_0{\rm span}\left\{\partial_zu_0,\partial_au_0,\partial_Eu_0\right\}=0,~~L_0\partial_cu_0=-\partial_zu_0
\end{equation}
where here the $T$-periodic boundary conditions have not yet been enforced; see Lemma \ref{l:gker} below for a precise description
of the generalized kernel of $L_0$.  

From the above remarks, it is natural that the spectral stability of the underlying wave $u_0$ to slow modulations 
corresponds to ensuring the spectrum of the Bloch operators $L_\xi$ near $\lambda=0$ is confined to the imaginary
axis for all $0\leq|\xi|\ll 1$.  Our program is thus two-fold: first study the structure
of the generalized kernel of the unmodulated operator $L_0$, and then use bifurcation theory to examine how
the spectrum near the origin of the modulated operator $L_\xi$ bifurcates from $\lambda=0$ for $0<|\xi|\ll 1$.  
This program is carried out in the next section.

\subsection{Modulational Instability in gKdV Type Equations}\label{s:migkdv}

We begin our program by demonstrating that, under a generic non-degeneracy
condition, the unmodulated operator $L_0$ has a two-dimensional $T$-periodic kernel with a Jordan
chain of height one.  Henceforth, we employ the notation
\[
\{f,g\}_{x,y}:=\det\left(\frac{\partial(f,g)}{\partial(x,y)}\right)\quad{\rm and}\quad\{f,g,h\}_{x,y,z}:=\det\left(\frac{\partial(f,g,h)}{\partial(x,y,z)}\right)
\]
for determinants of $2\times 2$ and $3\times 3$ Jacobians, respectively.

\begin{lemma}[\cite{Bronski.Johnson.2010,Bronski.Johnson.Kapitula.2011}]\label{l:gker}
Suppose that $u(x;a_0,E_0,c_0)$ is a $T$-periodic solution of the profile equation \eqref{quad2} and that the Jacobian determinants
$T_E$, $\{T,M\}_{a,E}$ and $\{T,M,P\}_{a,E,c}$ are non-zero at $(a_0,E_,c_0)$.  Then the functions
\[
\begin{array}{lcl}
\phi_0=\{T,u\}_{a,E} & \qquad\qquad &\psi_0=1\\
\phi_1=\{T,M\}_{a,E}u_x &\qquad\qquad & \psi_1=\int_0^x\phi_2(s)ds\\
\phi_2=\{u,T,M\}_{a,E,c} &\qquad\qquad &\psi_2=\{T,M\}_{a,E}+\{T,M\}_{a,E}u
\end{array}
\]
are all $T$-periodic and satisfy the relations
\[
\begin{array}{lcl}
L_0\phi_0=L_0\phi_1=0 & \qquad\qquad &  L_0^\dag\psi_0=L_0^\dag\psi_2=0\\
L_0\phi_2=-\phi_1 &  \qquad\qquad &  L_0^\dag\psi_1=\psi_2.
\end{array}
\]
In particular, the functions $\{\phi_j\}_{j=0}^2$ form a basis for the generalized null-space of $L_0$, and the functions
$\{\psi_j\}_{j=0}^2$ form a basis for the generalized null-space of $L_0^\dag$.  Moreover,
$\left<\psi_j,\phi_i\right>_{L^2_{\rm per}([0,T])}=0$ if $i\neq j$.
\end{lemma}

The proof follows from a straightforward application of the Fredholm alternative; see \cite{Bronski.Johnson.2010,Bronski.Johnson.Kapitula.2011}, for example.  The main observation is that while $u_x$ is $T$-periodic,
the derivatives $u_a$, $u_E$, and $u_c$ generically are not $T$-periodic due to the dependence of the period on the parameters
$(a,E,c)$. Indeed, the change in these functions across a period is proportional to derivatives of the period: for example,
\[
\left(\begin{array}{c}u_E(T)\\u_{Ex}(T)\\u_{Exx}(T)\\ \cdots\end{array}\right)-\left(\begin{array}{c}u_E(0)\\u_{Ex}(0)\\u_{Exx}(0)\\ \cdots\end{array}\right)
=T_E\left(\begin{array}{c}u_E(0)\\u_{Ex}(0)\\u_{Exx}(0)\\ \cdots\end{array}\right),
\]
with similar expressions for the change in $u_a$ and $u_c$ across a period.  It is then trivial to check that the functions $\{\phi_j\}_{j=0}^2$
are $T$-periodic and, thanks to \eqref{kerform}, satisfy the desired relations under the action of $L_0$.  Further, since $\phi_2$ clearly
has mean-zero, the functions $\{\psi_j\}_{j=0}^2$ are clearly $T$-periodic and, by again differentiating the profile equation, are seen
to satisfy the desired relations under the action of $L_0^\dag$.  Since $u$ can be chosen, through appropriate translation, to be an even function, and since
$L_0$ is parity reversing\footnote{That is, $L_0$ sends odd functions to even functions, and even functions to odd functions, all while preserving periodicity.},
we can further verify that $\left<\psi_j,L_0\phi_k\right>_{L^2_{rm per}([0,T])}=0$ for all $j,k=0,1,2$ which, by the Fredholm alternative,
clearly implies the desired structure of the generalized kernels of $L_0$ and $L_0^\dag$.

Henceforth, we assume that the quantities $T_E$, $\{T,M\}_{a,E}$, and $\{T,M,P\}_{a,E,c}$ are non-zero at the underlying $T$-periodic
traveling wave $u_0$.  Lemma \ref{l:gker} then implies that the elements of the $T$-periodic kernel of the unmodulated operator $L_0$
are given by the elements of the tangent space to the (two-dimensional) manifold of solutions \emph{of fixed period and fixed wave speed} at the point $u_0$,
while the generalized $T$-periodic kernel of $L_0$ consists also of a vector in the tangent space to the (three-dimensional) manifold of solutions
\emph{of fixed period} at $u_0$ with no restrictions on the wave speed.  It immediately follows that 
the origin is a $T$-periodic eigenvalue of $L_0$ with algebraic multiplicity three and geometric multiplicity two.

We now turn to study how this triple eigenvalues bifurcates from the $(\lambda,\xi)=(0,0)$ state.  To this end, we note that
a standard Baker-Campbell-Hausdorff expansion reveals that
\[
L_\xi=L_0+i\xi L_1+(i\xi)^2L_2+(i\xi)^3
\]
$L_0=\partial_x\mathcal{L}$ is the unmodulated operator, $L_1:=[L_0,x]=\mathcal{L}-2\partial_x^2$, and\footnote{Here and throughout,
$[A,B]=AB-BA$ denotes the commutator of operators $A$ and $B$.} $L_2:=[L_1,x]=-3\partial_x$.  Note
that these operators are well-defined on $L^2_{\rm per}([0,T])$ even though the function $x$ is not.
Furthermore, a very useful observation is that
\begin{equation}\label{property}
\left<\psi_i,L_j\phi_k\right>_{L^2_{\rm per}([0,T])}=0
\end{equation}
whenever $i+j+k=0$ mod$(2)$, which can readily be verified by using Lemma \ref{l:gker} and noting that $L_0$ and $L_2$ reverse parity
while $L_1$ preserves parity\footnote{That is, $L_0$ and $L_2$ send even (odd) functions to odd (even) functions, while $L_1$ sends
even (odd) functions to even (odd) functions.}.
Since $L_\xi$ is a relatively compact perturbation of $L_0$ depending analytically on the Bloch frequency $\xi$,
it follows that the operator $L_\xi$ will have three eigenvalues $\{\lambda_j(\xi)\}_{j=1}^3$, defined for $|\xi|\ll 1$,
bifurcating from $\lambda=0$  for $0<|\xi|\ll 1$.  To determine whether the underlying wave $u$ is modulationally stable or not,
it remains to track these eigenvalues as functions of $\xi$ and, in particular, determine whether these eigenvalues are confined to the
imaginary axis for all $0<|\xi|\ll 1$ or not.

Notice that by Lemma \ref{l:gker} we may use the functions $\{\phi_j\}_{j=0}^2$ and $\{\psi_j\}_{j=0}^2$ to construct explicit rank  3 eigenprojections
\[
\Pi_0:L^2_{\rm per}([0,T])\to{\rm gker}(L_0),~~~\widetilde{\Pi}_0:L^2_{\rm per}([0,T])\to{\rm gker}(L_0^\dag)
\]
onto the total eigenspace\footnote{Here, the total eigenspaces for $L_0$ and $L_0^\dag$ are precisely ${\rm gker}(L_0)={\rm ker}(L_0^2)$
and ${\rm gker}={\rm ker}(L_0^2)$, respectively.} for the operators $L_0$ and $L_0^\dag$.
We can thus represent the action of the operator $L_0$ on generalized kernel of $L_0$ as the $3\times 3$ matrix operator\footnote{Specifically,
$M_0:=\left[\langle\psi_i,L_0\phi_j\rangle>_{L^2_{\rm per}([0,T])}\right]_{i,j=1,2,3}$ and similarly for the projection of the identity; see below.  In Section \ref{s:wsmall} we will
utilize a slightly different rescaling of these matrix operators.}
\[
M_0:=\widetilde\Pi_0 L_0\Pi_0
=\left(\begin{array}{ccc}
			0 & 0 & 0\\
			0 & 0 & \left<\psi_1,L_0\phi_2\right>\\
			0 & 0 & 0
			\end{array}\right).
\]
Notice that $\left<\psi_1,L_0\phi_2\right>=\frac{1}{2}\{T,M\}_{a,E}\{T,M,P\}_{a,E,c}\neq 0$ by assumption, reflecting the Jordan structure of the unmodulated
operator $L_0$.  Further, the above projections are non-degenerate in the sense the action of the identity operator
\[
\widetilde\Pi_0\Pi_0=\left(\begin{array}{ccc}
											\left<\psi_0,\phi_0\right> & 0 & 0\\
											0 & \left<\psi_1,\phi_1\right> & 0\\
											0 & 0 & \left<\psi_2,\phi_2\right>
											\end{array}\right),
\]
is non-singular with determinant $-\{T,M\}_{a,E}^3\{T,M,P\}_{a,E,c}$.

Since  $\lambda=0$ is an \emph{isolated} eigenvalue of $L_0$, it follows by standard spectral perturbation theory that the above eigenprojections
can be continued for $|\xi|\ll 1$ into analytically varying rank 3 projections
\[
\Pi_\xi:L^2_{\rm per}([0,T])\to\bigoplus_{j=1}^3{\rm ker}\left(L_\xi-\lambda_j(\xi){\bf I}\right),~~
\widetilde{\Pi}_\xi:L^2_{\rm per}([0,T])\to\bigoplus_{j=1}^3{\rm ker}\left(L_\xi^\dag-\overline{\lambda_j(\xi)}{\bf I}\right)
\]
with ranges coinciding with the total left and right eigenspaces associated  with the eigenvalues $\{\lambda_j(\xi)\}_{j=1}^3$ of $L_\xi$; see \cite{Kato.book}.  
In particular, we may find analytically 
varying bases $\{v_j(\xi)\}_{j=0}^2$ and $\{\tilde{v}_j(\xi)\}_{j=0}^2$ of the total right and left eigenspaces for the eigenvalues  $\{\lambda_j(\xi)\}_{j=1}^3$
in such a way that $v_j(0)=\phi_j$ and $\tilde{v}_j(0)=\psi_j$.  Note in particular that we need not require the $\{v_j\}_{j=0}^2$ and $\{\tilde{v}_j\}_{j=0}^2$
to be \emph{eigen-bases} for the total eigenspace, a degree of flexibility that is helpful during explicit calculations.  
Using these $\xi$-dependent projections, for each fixed $|\xi|\ll 1$ 
we can project the infinite dimensional spectral problem \eqref{spec} onto the three dimensional total eigenspace
associated with the  $\{\lambda_j(\xi)\}_{j=1}^3$ and, in particular, the action of the operators $L_\xi$ on this subspace
can be represented by the $3\times 3$ matrix operator\footnote{As above, $M_\xi:=\left[\langle\tilde{v}_j(\xi),L_\xi v_j(\xi)\rangle_{L^2_{\rm per}([0,T])}\right]_{i,j=1,2,3}$ and similarly
with the projection of the identity operator.}
\[
M_\xi:=\widetilde{\Pi}_\xi L_\xi\Pi_\xi.
\]
The key observation is that for each fixed $\xi$, $M_\xi$ is a constant matrix and that, further, 
the eigenvalues of the matrix $M_\xi$ coincide \emph{precisely} with the eigenvalues  $\{\lambda_j(\xi)\}_{j=1}^3$
for the modulated operators $L_\xi$.  In particular, for each $|\xi|\ll 1$ the eigenvalues $\{\lambda_j(\xi)\}_{j=1}^3$ are the roots of the cubic
characteristic equation
\[
\det\left(M_\xi-\lambda\widetilde{\Pi}_\xi\Pi_\xi\right)=0.
\]
Note that since $\widetilde\Pi_0\Pi_0$ is invertible, it follows by continuity
that $\widetilde\Pi_\xi\Pi_\xi$ is non-singular for all $|\xi|\ll 1$ .
Thus, the three eigenvalues bifurcating from the triple eigenavlue $\lambda=0$ of $L_0$  can be studied for for $|\xi|\ll 1$
by analyzing the eigenvalues of the $3\times 3$ matrix $M_\xi$, which is a seemingly much easier task.

The eigenvalues of the matrix $M_\xi$ may now be studied using standard techniques in matrix perturbation theory.  First,
the matrices $M_\xi$ depend analytically on $\xi$ for $|\xi|\ll 1$, due to the analyticity of the associated eigenprojections $\Pi_\xi$ and $\widetilde{\Pi}_\xi$
and the analytic dependence of the operators $L_\xi$ on $\xi$, and can be expanded for $|\xi|\ll 1$ as
\[
M_\xi = M_0 + (i\xi) M^{(1)} + (i\xi)^2M^{(2)}+\mathcal{O}(|\xi|^3).
\]
for some $3\times 3$ matrices $M^{(1)}$ and $M^{(2)}$ whose coefficients are independent of $\xi$.  It is well known that the eigenvalues of $M_\xi$ in general
bifurcate continuously in the perturbation parameter $\xi$, but may \emph{not} vary in a $C^1$ manner.  Rather, due to the Jordan block
at $\xi=0$, these eigenvalues are generically expected to 
admit a Puiseaux series expansion in fractional powers of $\xi$.  In order to guarantee that the eigenvalues are at least $C^1$ in $\xi$ near $\xi=0$,
it is sufficient to verify that the 1,2 and 3,2 entries on the matrix $M^{(1)}$ both vanish.  By a straightforward calculation, it can be shown that
\[
[M^{(1)}]_{1,2}=\left<\psi_0,L_1\phi_1\right>_{L^2_{\rm per}([0,T])}\quad\textrm{and}\quad [M^{(1)}]_{3,2}=\left<\psi_2,L_1\phi_1\right>_{L^2_{\rm per}([0,T])},
\]
both of which vanish thanks to \eqref{property}.  Consequently, the eigenvalues $\{\lambda_j(\xi)\}_{j=1}^3$ are $C^1$ functions of $\xi$ and hence
admit an expansion of the form
\begin{equation}\label{evexpand}
\lambda_j(\xi)=i\xi\mu_j(\xi)=i\xi\mu_j(0)+o(\xi),~~~|\xi|\ll 1,
\end{equation}
for some $C^1$ function $\mu_j(\xi)$ determined by the roots of the characteristic equation
\begin{equation}\label{charinform}
\det\left(\frac{1}{i\xi}M_\xi-\mu \widetilde\Pi_\xi\Pi_\xi\right)=0
\end{equation}
for $|\xi|\ll 1$.
It follows that, to leading order, the modulational stability of the underlying wave may be inferred from
determining the values of $\mu_j(0)\in\mathbb{C}$ in \eqref{evexpand}.
Note that if we further assume the numbers $\mu_j(0)$ are distinct (which generically holds in our case), then the eigenvalues
$\{\lambda_j(\xi)\}_{j=1}^3$ depend analytically on $\xi$ for $|\xi|\ll 1$.

To determine the values $\mu_j(0)$ in \eqref{evexpand}, notice the Dunford integral calculus implies that
the bases $\{v_j(\xi)\}_{j=1}^3$ and $\{\tilde v_j\}_{j=0}^2$ have the same regularity on $\xi$ as the eigenvalues $\{\mu_j\}_{j=1}^3$,
and hence they may be expanded near  $\xi=0$ as
\begin{equation}\label{basisexpand}
v_j(\xi)=\phi_j+(i\xi)v_j^{(1)}+o(\xi),~~\tilde v_j(\xi)=\psi_j+(i\xi)\tilde v_j^{(1)}+o(\xi)
\end{equation}
for some functions $v_j^{(1)}$ and $\tilde v_j^{(1)}$ independent of $\xi$.  We wish to use these expansions to compute
asymptotic expansions of the projections $\Pi_\xi$ and $\widetilde{\Pi}_\xi$ valid near $\xi=0$.  Note, however, the non-trivial
Jordan structure of $M_0$ implies an inherent degeneracy in the limit $\xi\to 0$.  Indeed,
note that while the kernel of $M_0$ is two-dimensional, there are three eigenvalues of $M_\xi$ defined near the origin
for $\xi=0$.  Moreover, the three eigenvectors of $M_\xi$ tend to \emph{the same limit}\footnote{This is a consequence
of the fact that all three of the eigenfunctions for $L_\xi$ associated to the eigenvalues $\{\lambda_j(\xi)\}_{j=1}^3$ 
tend to $\phi_1$ as $\xi\to 0$, a fact seen both analytically and numerically.} as $\xi\to 0$.
An effecient way of unfolding this degeneracy is to rescale the matrices $(i\xi)^{-1}M_\xi$ and $\widetilde\Pi_\xi\Pi_\xi$ as
\[
\widehat{M}_{\xi}:=\frac{1}{i\xi}\Sigma(\xi)^{-1}M_\xi\Sigma(\xi),\quad \widehat{I}_\xi:=\Sigma(\xi)^{-1}\widetilde\Pi_\xi\Pi_\xi\Sigma(\xi),
\]
where
\[
\Sigma(\xi):=\left(\begin{array}{ccc}
								i\xi & 0 & 0\\
								0 & 1 & 0\\
								0 & 0 & i\xi
								\end{array}\right).
\]
By straightforward computations, it follows that for $|\xi|\ll 1$ we have the expansions
\[
\widehat{M}_\xi=\left(\begin{array}{ccc}
			\left<\psi_0,L_1\phi_0\right> & m_{1,2} & 
																											\left<\psi_0,L_1\phi_2\right>+\left<\tilde v_0^{(1)},L_0\phi_2\right>\\
			0 & \left<\psi_1,L_1\phi_1+L_0v_1^{(1)}\right> & \left<\psi_1,L_0\phi_2\right>\\
			\left<\psi_2,L_1\phi_0\right> & m_{3,2} &
																											\left<\psi_2,L_1\phi_2\right>+\left<\tilde v_2^{(1)},L_0\phi_2\right>
							\end{array}\right)+o(1)
\]
where
\[
m_{j,2}:=\left<\psi_j,L_2\phi_1+L_1v_1^{(1)}\right>+\left<\tilde v_j^{(1)},L_1\phi_1+L_0v_1^{(1)}\right>,~~j=1,3
\]
and that, similarly,
\[
\widehat{I}_\xi=
\left(\begin{array}{ccc}
\left<\psi_0,\phi_0\right> & \left<\psi_0,v_1^{(1)}\right>+\left<\tilde v_0^{(1)},\phi_1\right> & 0\\
0 & \left<\psi_1,\phi_1\right> & 0\\
0 & \left<\psi_2,v_1^{(1)}\right>+\left<\tilde v_2^{(1)},\phi_1\right> & \left<\psi_2,\phi_2\right>
\end{array}\right) + o(1).
\]
Therefore, it follows that the numbers $\mu_j(0)$ are determined as the roots of the (cubic) \emph{effective dispersion relation}
\[
\mathcal{D}(\mu):=\det\left(\widehat{M}_0-\mu \widehat{I}_0\right)=0
\]
This establishes the following theorem, which is a restatement of Theorem 2 in \cite{BrJZ}.

\begin{theorem}[\cite{BrJZ}]\label{t:whitham}
Let $u(\cdot;a_0,E_0,c_0)$ be a $T=T(a_0,E_0,c_0)$ periodic traveling wave solution of \eqref{gkdv} and suppose that
$T_E$, $\{T,M\}_{a,E}$, and $\{T,M,P\}_{a,E,c}$ 
are non-zero at $a_0$, $E_0$, and $c_0$.
Then for $|\xi|\ll 1$, the triple eigenvalue $\lambda=0$ of $L_0$ bifurcates
into $3$ eigenvalues $\lambda_j(\xi)$ that are $C^1$ in $\xi$ near $\xi=0$ and satisfy
\[
\lambda_j(\xi)=i\mu_j\xi+o(\xi),~~j=1,2,3
\]
where the constants $\mu_j\in\mathbb{C}$ are the eigenvalues of the effective dispersion matrix 
\begin{equation}\label{edm}
{\bf D}(a,E,c):=\left(\left(\widehat{I}_0\right)^{-1}\widehat{M}_0\right)(a,E,c)
\end{equation}
at $(a,E,c)=(a_0,E_0,c_0)$.
In particular, a necessary condition for $u$ to be modulationally stable is that the dispersion matrix
be weakly hyperbolic, i.e. has three real eigenvalues.
\end{theorem}

To utilize Theorem \ref{t:whitham} in practice, one must determine the first-order correctors $\tilde{v}_0^{(1)}$, $\tilde{v}_2^{(1)}$, and $v_1^{(1)}$
in \eqref{basisexpand}.  We note, however, that in the \emph{small amplitude limit} these correctors do not contribute
to leading order and hence are not necessary to compute the roots of $\mathcal{D}$; see Section \ref{s:wsmall} below.  
Outside of this asymptotic regime, however, these correctors play an important role.  First, notice that
the  functions $\{v_j(\xi)\}_{j=0}^2$ and $\{\tilde v_j(\xi)\}_{j=0}^2$  satisfy the relations
\[
\Pi_\xi\left(L_\xi v_j(\xi)\right)=L_\xi v_j(\xi),\quad \widetilde{\Pi}_\xi\left(L_\xi^\dag \tilde v_j(\xi)\right)=L_\xi^\dag\tilde v_j(\xi),
\]
respectively, for all $|\xi|\ll 1$.  In particular, expanding in $\xi$ we find that
\[
L_0v_j^{(1)}+\left(L_1+\Pi_1L_0\right)\phi_j\in{\rm Range}\left(\Pi(0)\right)
\]
and
\[
L_0^\dag\tilde{v}_j^{(1)}+\left(L_1^{\dag}+\widetilde\Pi_1L_0^\dag\right)\psi_j\in{\rm Range}\left(\widetilde\Pi(0)\right)
\]
In particular, since $L_0\phi_j=0$ we may choose $v_1^{(1)}$ so that $L_0v_1^{(1)}=-L_1\phi_1$.
An explicit formula for $v_1^{(1)}$ can be found by noting that, by definition,
\begin{equation}\label{Po}
L_1u_x=[L_0,x]u_x=L_0\left(xu_x\right)
\end{equation}
and that $T_E xu_x+Tu_E$ is $T$-periodic.
Similarly, since ${\rm ker}(L_0^\dag)={\rm span}\{\psi_j\}_{j=0,2}$ we may choose $\tilde{v}_0^{(1)}$ and $\tilde{v}_2^{(1)}$ such that
$L_0^\dag\tilde{v}_j^{(1)}=-L_1^\dag\psi_j$ for $j=0,2$.  Unlike $v_1^{(1)}$, explicit forms of $\tilde{v}_j^{(1)}$, $j=0,2$
are not needed; their defining relations are sufficient for calculation purposes.
With these choices one can, through tedious but straightforward calculations, compute the effective dispersion relation \emph{explicitly}, yielding
\[
\mathcal{D}(\mu)=
C\left(-\mu^3+\frac{\mu}{2}\left(\{T,P\}_{E,c}+2\{M,P\}_{a,E}\right)-\frac{1}{2}\{T,M,P\}_{a,E,c}\right)
\]
for some non-zero constant $C=C(a,E,c)$.  The reality of the roots of $\mathcal{D}$ follows from the sign of its discriminant, given here
as
\[
\Delta_{MI}:=\frac{1}{2}\left(\{T,P\}_{E,c}+2\{M,P\}_{a,E}\right)^3-\frac{27}{4}\{T,M,P\}_{a,E,c}^2.
\]
In particular, the polynomial $\mathcal{D}$ has three distinct real roots provided that $\Delta_{MI}>0$, while it has one real and two (non-real) complex
conjugate roots if $\Delta_{MI}<0$.  It follows that if $\Delta_{MI}<0$, then $u_0$ is modulationally unstable with one branch of spectrum bifurcating
from the origin along the imaginary axis two branches bifurcating in complex directions.  
It follows that a necessary condition for modulational stability is $\Delta_{MI}\geq 0$.  Further, when $\Delta_{MI}>0$,
a simple symmetry argument implies that the three spectral branches bifurcating from the origin are confined (locally)
to the imaginary axis, corresponding to a triple covering, yielding modulational stability.

\begin{corollary}[\cite{Bronski.Johnson.2010}]\label{c:MI}
Under the hypotheses of Theorem \ref{t:whitham}, a necessary condition for the modulational of the underlying wave $u$
$\Delta_{MI}\geq 0$.   Furthermore, a sufficient condition for $u$ to be modulationally stable is that $\Delta_{MI}>0$.
\end{corollary}

\subsection{Connection to Whitham Modulation Theory}\label{s:validate}

In Section \ref{s:whitham} we used formal asymptotics to determine a dispersion relation $\widehat{\mathcal{D}}_W$ in \eqref{Wdisp} with the property
the linearized Whitham system \eqref{linwhitham} has eigenvalues of the form $\lambda=ik \mu$ whenever $\mu\in\mathbb{C}$
satisfies $\widehat{\mathcal{D}}_W(z)=0$.  In Section \ref{s:migkdv} above, we demonstrated that he eigenvalues for the Bloch operators $L_\xi$ expand as $\lambda=i\mu\xi+o(\xi)$ where $\mu\in\mathbb{C}$
is one of the three (distinct by assumption) eigenvalues of the effective dispersion matrix ${\bf D}(a,E,c)$ defined in \eqref{edm}.
These two approaches are in fact equivalent for equations of gKdV type.

\begin{theorem}[\cite{JZ,BrJZ}]\label{t:validate}
Under the assumptions of Theorem \ref{t:whitham}, the eigenvalues of the effective dispersion matrix \eqref{edm}
agree with the roots of the reduced Whitham dispersion relation $\widehat{D}_W(z)=0$.  In particular, a necessary condition
for the modulational stability of the underlying periodic wave train is that the linearized Whitham system \eqref{linwhitham} be
weakly hyperbolic, or, equivalently, that equation $\widehat{\mathcal{D}}_W(z)=0$ has three real roots.  Further, a sufficient
condition for modulational stability is for \eqref{linwhitham} to be strictly hyperbolic, i.e. the equation $\widehat{\mathcal{D}}_W(z)=0$
has three distinct real roots.
\end{theorem}

Theorem \ref{t:validate} rigorously justifies formal predictions of Whitham's theory of modulations in the context of  gKdV equations
of the form \eqref{gkdv}.  It follows that a sufficient condition for the modulational instability of the unerlying carrier wave is the \emph{ellipticity} of the
associated linearized Whitham system.  Its proof follows by a long and tedious calculation demonstrating that there exists a constant $C\neq 0$ such that
\[
\det\left({\bf D}(a,E,c)-\mu{\bf I}\right)=C\widehat{\mathcal{D}}_W(\mu)
\]
for all $\mu\in\mathbb{C}$: see \cite{JZ,BrJZ} for details.  We merely pause to note that, in a sense, 
Theorem \ref{t:validate} implies that the processes of averaging and linearizing commute: the linearized Whitham system being obtained from
and averaging process (Fredholm alternative) and then linearizing, while the rigorous theory obtained from linearizing the PDE and then
applying the Fredholm alternative (averaging).  This observation suggests the above program can be readily extended
to equations outside the context of \eqref{gkdv}

One class of equations for which the above program has been particularly successful is in the context of dissipative PDE of the form
\[
u_t+f(u)_x+g(u)=\left(B(u)u_x\right)_x,~~x\in\mathbb{R},~~u\in\mathbb{R}^n
\]
consisting of a continuum of models, with $f=0$ corresponding to systems of reaction diffusion equations, and $g=0$ to systems of conservation laws.
Not only has Whitham's theory of modulations been rigorously validated in this context at the level of spectral stability, it has also provided key insights
into the \emph{nonlinear dynamics} of modulated waves in such models, as well as the stability and dynamics of wave trains to \emph{nonlocalized perturbations},
consisting of a localized part plus an asymptotic phase shift at the spatial infinities; see \cite{Johnson.Noble.Rodrigues.Zumbrun2015,BJNRZ.2013,Johnson.Noble.Rodrigues.Zumbrun.2014,Johnson.Noble.Rodrigues.Zumbrun.2013a,Johnson.Noble.Rodrigues.Zumbrun.2013b}.  These additional insights come from
continuing the asymptotic WKB expansion in Section \ref{s:whitham} to next order in $\varepsilon$, obtaining a diffusive correction to the Whitham
modulation system.  The details of this program are beyond the scope of the current manuscript, and the reader is referred to the above references for details.

\subsection{Evaluation of $\Delta_{MI}$}

In Section 4.2 above, we defined a modulational instability index $\Delta_{MI}$,  the sign of which determines the modulational stability of a 
given periodic traveling wave solution of the gKdV equation \eqref{gkdv}.  To compute $\Delta_{MI}$ in practice, we must
determine derivatives of the quantities $T$, $M$, and  $P$ with respect to the traveling wave parameters $a$, $E$, and $c$,
a task that in general seems formidable.  Nevertheless, when $f$ is a polynomial nonlinearity these derivatives may be 
computed \emph{explicitly} in terms of the underlying wave.
Indeed, note that if $f$ is a polynomial then the quantities $T$, $M$, and $P$ can be expressed in terms of the moments
\[
\zeta_k:=\int_0^Tu^k(x)dx=\oint_\gamma\frac{u^k}{\sqrt{E-V(u;a,c)}}~du,
\]
where here $E-V(u;a,c)$ is a polynomial of one degree higher than $f$; indeed, $T=\gamma_0$, $M=\gamma_1$, and $P=\gamma_2$.  
Further, the derivatives of $T$, $M$, and $P$ may be expressed in terms of the moments
\[
I_k:=\oint_\gamma\frac{u^k}{\left(E-V(u;a,c)\right)^{3/2}}~du;
\]
for example, $T_E=I_0$ and $P_c=\frac{1}{2}I_2$.  The amazing fact is that if $f$ is a polynomial of degree $n$,
then the $2n-1$ moments $\{I_k\}_{k=0}^{2n-2}$ can be expressed as linear combinations of the $n-1$ moments
$\{\zeta_k\}_{k=0}^{n-2}$, a fact known as the Picard-Fuchs relation.  As such, the modulational
instability index $\Delta_{MI}$ can be \emph{explicitly computed} in terms of moments of the underlying wave $u$;
see \cite{Bronski.Johnson.Kapitula.2011} for more details.

To illustrate how this process works, suppose $f$ is a polynomial of degree $n-1$ and set $P(u)=E-V(u;a,c)$.
Then we can write $P(u)=a_0+a_1u+a_2u^2+\ldots a_nu^n$ for some constants $\{a_j\}_{j=0}^n$.  Defining
the moments $\zeta_k$ and $I_k$ as above, it clearly follows that
\[
\frac{d\zeta_k}{da_j}=\frac{d\zeta_j}{da_k}=-\frac{1}{2}I_{k+j}
\]
for each $k,j=0,1,2,\ldots,n$.  Upon observing that
\[
\zeta_k=\oint_\gamma\frac{u^kP(u)}{P(u)^{3/2}}~du=\sum_{j=0}^na_jI_{j+m},~~m=0,1,2,\ldots,n-1
\]
and that\footnote{Here, we take the convention that the left hand side is zero when $k=0$.}
\[
2k\zeta_{k-1}=2k\oint_\gamma\frac{u^{k-1}}{\sqrt{P(u)}}~du=\oint_\gamma\frac{u^kP'(u)}{P(u)^{3/2}}~du
=\sum_{j=0}^nja_jI_{j+k-1},~~k=0,\ldots,n-1
\]
we find the following linear system of $2n-1$ equations in the unknowns $\{I_k\}_{k=0}^{2n-2}$:
\begin{equation}\label{PF}
\left(\begin{array}{ccccccc}
a_0 & a_1 & \cdots & a_n & 0 & 0 & \cdots\\
0 & a_0 & a_1 & \cdots & a_n & 0 & \cdots\\
\vdots & \ddots & \ddots & \ddots & \ddots & \ddots & \ddots\\
0 & \cdots & 0 & a_0 & a_1 & \cdots & a_n\\
a_1 & 2a_2 & \cdots & na_n & 0 & 0 & \cdots\\
0 & a_1 & 2a_2 & \cdots & na_n & 0 & \cdots\\
\vdots & \ddots & \ddots & \ddots & \ddots & \ddots & \ddots\\
0 & \cdots & 0 & a_1 & 2a_2 & \cdots & na_n
\end{array}\right)
\left(\begin{array}{c} I_0 \\ I_1 \\ I_2 \\ I_3 \\ \vdots \\ \vdots \\ I_{2n-3} \\ I_{2n-2}\end{array}\right)
=
\left(\begin{array}{c} \zeta_0 \\ \zeta_1 \\ \vdots \\ \zeta_{n-2} \\ 0 \\ 2\zeta_0 \\ \vdots \\ 2(n-1)\zeta_{n-2}\end{array}\right),
\end{equation}
which we refer to as the Picard-Fuchs system.
The matrix that arises in the above linear system may be recognized as the Sylvester matrix of the polynomials $P(u)$ and $P'(u)$,
A standard result in commutative algebra states that the Sylvester matrix of two polynomials is non-singular if and only if the two polynomials
have no common roots.  Thus, the above matrix is invertible provided that the polynomial $P$ has only simple roots, in which case
the moments $\{I_k\}_{k=0}^{2n-2}$ can be explicitly expressed as a linear combination of the $\{\zeta_k\}_{k=0}^{n-2}$ with coefficients
explicitly determined by the coefficients of $P$.

For a given polynomial $P$, the above linear system can be solved easily with the aid of a computer algebra system.  
In the next section, we discussion applications of the above procedure in the case where $P(u)=E-V(u;a,c)$ and $V$
is the effective potential energy associated to \eqref{gkdv} for various nonlinearities $f$.

\section{Applications}\label{s:app}

In Section 4 above, we defined a modulational instability index $\Delta_{MI}$,  the sign of which determines the modulational stability of a 
given periodic traveling wave solution of the gKdV equation \eqref{gkdv}.  We now discuss a few examples where one can
use the Picard-Fuchs system \eqref{PF} to explicitly evaluate $\Delta_{MI}$ in terms of the underlying wave.  We note that
these examples are \emph{not exhaustive} but are chosen for illustrative  purposes.

\subsection{The KdV Equation}

Consider first the KdV equation
\begin{equation}\label{kdv}
u_t=u_{xxx}+(u^2)_x.
\end{equation}
In this case, using notation from Section \ref{s:exist}, the effective potential energy is
\[
V(u;a,c)=-au-\frac{c}{2}u^2+\frac{1}{3}u^3
\]
and the periodic traveling waves $u(x;a,E,c)$ of \eqref{kdv} 
exist provided that the discriminant
\[
{\rm disc}\left(E-V(\cdot;a,c)\right)=\frac{1}{12}\left(16a^3+3a^2c^2-36Eac-6Ec^3-36E^2\right)
\]
is positive\footnote{Positivity of the discriminant ensures the polynomial $E-V(u;a,c)$ has three real roots in $u$.}.  Further, the only periodic solutions 
along the surface ${\rm disc}\left(E-V(\cdot;a,c)\right)=0$ are constants.
In this case, setting
\[
\zeta_k:=\oint_\gamma\frac{u^k~du}{\sqrt{2(E-V(u;a,c)}},\quad I_k=\oint_\gamma\frac{u^k~du}{\left(2\left(E-V(u;a,c)\right)\right)^{3/2}},
\]
we find that $(T,M,P)=(\zeta_0,\zeta_1,\zeta_2)$ and that
\[
\nabla_{E,a,c}\zeta_k=\left<-I_k,-I_{k+1},-\frac{1}{2}I_{k+2}\right>.
\]
In particular, the Picard-Fuchs system \eqref{PF} becomes
\[
\left(\begin{array}{ccccc}
E & a & \frac{c}{2} & \frac{1}{3} & 0\\
0 & E & a & \frac{c}{2} & \frac{1}{3}\\
a & c & -1 & 0 & 0\\
0 & a & c & -1 & 0\\
0 & 0 & a & c & -1
\end{array}\right)
\left(\begin{array}{c}I_0\\ I_1\\ I_2\\ I_3\\ I_4\end{array}\right)
=
\left(\begin{array}{c}T \\ M \\0 \\ 2T \\ 4M\end{array}\right)
\]
which, after some elementary matrix algebra, implies the following  formulae for 
\begin{align*}
&T_E=\frac{(4a+c^2)M+(6E+ac)T}{12{\rm disc}\left(E-V(\cdot;a,c)\right)},\qquad
\{T,M\}_{a,E}=-\frac{T^2V'(M/T;a,c)}{12{\rm disc}\left(E-V(\cdot;a,c)\right)}\\
&\qquad\qquad\qquad\{T,M,P\}_{a,E,c}=\frac{T^3(E-V(M/T;a,c))}{2{\rm disc}\left(E-V(\cdot;a,c)\right)}\\
\end{align*}
and
\[
2\Delta_{MI}=\frac{\left(\alpha_{3,0}T^3+\alpha_{2,1}T^2M+\alpha_{1,2}TM^2+\alpha_{0,3}M^3\right)^2}{2^{10}3^7{\rm disc}\left(E-V(\cdot;a,c)\right)^3};
\]
see \cite{Bronski.Johnson.Kapitula.2011}.
Here the coefficients $\alpha_{ij}$ are given by 
\[
\begin{array}{lcl}
\alpha_{3,0}= 36 E + 18 a E c - 8 a^3, &\quad& \alpha_{2,1}= 18 E c^2 - 6 a^2 c + 36 a E, \\
\alpha_{2,1}= -18 c E + 24 a^2 + 3 a c^2, &\quad& \alpha_{0,3}= c^3 + 6 a c + 12 E.
\end{array}
\]
The monotonicity of the period, $T_E>0$, holds by a result of Schaaf\cite{Schaaf.1985}, while $\{T,M,P\}_{a,E,c}>0$ clearly holds.  Further, $\{T,M\}_{a,E}$ is positive due to Jensen's inequality
and the fact that
\[
\oint\frac{V'(u;a,c)~du}{\sqrt{E-V(u;a,c)}}=0.
\]
Thus, the hypothesis of Theorem \ref{t:whitham} are satisfied, and hence Corollary \ref{c:MI} implies that the sign of $\Delta_{MI}$
determines the modulational stability of any underlying wave $u$.  Since the numerator of $\Delta_{MI}$ is seen numerically not to vanish,
it follows that $\Delta_{MI}>0$.  Taken together, this yields the following.

\begin{theorem}[\cite{Bronski.Johnson.Kapitula.2011}]
All periodic traveling wave solutions of the KdV equation \eqref{kdv} are modulationally stable.
\end{theorem}

\subsection{The Modified KdV Equation}

Next, we consider the modified KdV (mKdV) equation
\begin{equation}\label{mkdv}
u_t=u_{xxx}\pm(u^3)_x,
\end{equation}
which corresponds to the focusing (+) and defocusing (-) modified KdV (mKdV) equations.
Like the KdV equation, \eqref{mkdv} is completely integrable and hence its sepctrum can, in principle, be explicitly computed with the aid
of completely integrability and algebro-geometric techniques.
Of course, as mentioned above, our methods are independent of integrability.  Further, as we will see
our methods yield a modulational instability result for \eqref{mkdv} that seems to have been previously missed.

As discussed in Section \ref{s:exist} above, the existence of periodic solutions for \eqref{mkdv} differs between the focusing and defocusing cases.
In each case, however, existence may be inferred from properties of the associated effective potential energy
\[
V_{\pm}(u;a,c):=\mp \frac{1}{4}u^4-\frac{c}{2}u^2-au
\]
Indeed, in the focusing case, periodic solutions exist provided that the quartic polynomial $E-V_+(\cdot;a,c)$ has at least two real roots
and if its discriminant is non-zero.  In the defocusing case, however, periodic solutions only exist when the discriminant
of $E-V_-(\cdot;a,c)$ is positive, indicative of the polynomial having four distinct real roots.  See Figure \ref{fig1} in Section \ref{s:exist}.

In this case, the Picard-Fuchs system reads
\[
\left(\begin{array}{ccccccc}
E & a & c/2 & 0 & \mp 1/4 & 0 & 0\\
0 &E & a & c/2 & 0 & \mp 1/4 & 0\\
0 & 0 & E & a & c/2 & 0 & \mp 1/4 \\
a & c & 0 & \mp 1 & 0 & 0 & 0\\
0 & a & c & 0 & \mp 1 & 0 & 0\\
0 & 0 & a & c & 0 & \mp 1 & 0\\
0 & 0 & 0 & a & c & 0 & \mp1
\end{array}\right)
\left(\begin{array}{c} I_0 \\ I_1 \\ I_2 \\ I_3 \\ I_4 \\ I_5 \\ I_6\end{array}\right)
=
\left(\begin{array}{c} T \\ M \\ P \\ 0 \\ 2T \\ 4M \\ 6P\end{array}\right),
\]
which can be readily solved to provide explicit formulas for $T_E$, $\{T,M\}_{a,E}$, and $\{T,M,P\}_{a,E,c}$
in terms of $T$ and $P$ and\footnote{The fact that $M$ does \emph{not} enter the formulas for the gradients
of $T,M,P$ with respect to $a,E,c$ is easily explained by complex analysis; see \cite{Bronski.Johnson.Kapitula.2011} for details.} the parameters $(a,E,c)$; see \cite{Bronski.Johnson.Kapitula.2011}.
Further, we find
\[
\Delta_{MI}=\frac{\Gamma_{mKdV}(T,P,a,E,c)^2}{{\rm disc}\left(E-V_{\pm}(\cdot;a,c)\right)^3}
\]
where  $\Gamma_{mKdV}$is a homogeneous polynomial of degree 6 in $T$ and $P$ with coefficients depending on $a$, $E$, and $c$.
The sign of this clearly agrees with the sign of the discriminant of the quartic polynomial $E-V_{\pm}(\cdot;a,c)$, which in turn is positive
if it has four distinct real roots, and negative if it has only two distinct real roots and two distinct complex conjugate roots.
This leads us to the surprisingly simple characterization of modulational instability for the mKdV equation.

\begin{theorem}[\cite{Bronski.Johnson.Kapitula.2011}]\label{stability:mkdv}
The traveling wave solutions to the mKdV equation \eqref{mkdv} are modulationally unstable for a given set of parameters
$(a,E,c)$  if the polynomial
\[
\mathcal{P}(u;a,E,c)=E+au+\frac{c}{2}u^2\pm\frac{1}{4}u^4
\]
has only two distinct real roots and two distinct complex conjugate roots, and it is modulationally stable if it has four distinct real roots.
\end{theorem}

As mentioned above, the structure of the roots of $\mathcal{P}(\cdot;a,E,c)$ depends greatly on whether the equation is focusing or defocusing.
In the case of the focusing mKdV, Theorem \ref{stability:mkdv} implies that if the parameter values give rise to only one periodic solution,
then this solution is modulationally unstable.  Conversely, if the parameters give rise to two periodic solutions, then both of the solutions
are modulationally stable.  In terms of the explicit solutions described in Section \ref{s:explicit} above, we have the following corollary.

\begin{corollary}
For the focusing mKdV equation, all the dnoidal wave solutions, constructed in Section \ref{s:explicit} above, are modulationally stable, while all the cnoidal waves
are modulationally unstable.  
\end{corollary}

In contrast, notice that for periodic waves to exist in the defocusing mKdV
$\mathcal{P}$ \emph{must have four distinct zeroes}. This simple observation yields the following.

\begin{corollary}\label{c:dfmkdv}
All periodic traveling wave solutions of the defocusing mKdV equation are modulationally stable.  
\end{corollary}

In particular, Corollary \ref{c:dfmkdv} implies the snoidal wave solutions described in Section \ref{s:explicit} above are modulationally stable.

\subsection{The Schamel Equation}
Lastly, we study the Schamel equation (after taking $u\mapsto -u$ in \eqref{kdv})
\begin{equation}\label{kdv1/2}
u_t+u_{xxx}+\frac{5}{2}\left(|u|^{3/2}\right)_x=0,
\end{equation}
where $u$ must be taken to be positive. 
Periodic traveling waves of \eqref{kdv1/2} are determined implicitly via
\[
\oint\frac{du}{\sqrt{2(E+au+\frac{1}{2}cu^2-u^{5/2})}}=x-ct.
\]
Using the obvious change of variables $u=v^2$ shows that the traveling wave solution to this equation
may in turn be found from 
\[
\frac{v_x^2}{2}=E+av^2+\frac{c}{2}v^4-v^5.
\]
Thus, the physically admissible parameter regime corresponds to the set of all $(a,E,c)\in\mathbb{R}^3$ for which
there exists a bounded interval in $(0,\infty)$ such that the quintic polynomial $E+av^2+\frac{c}{2}v^4-v^5$ 
is non-negative, corresponding to a non-negative periodic solution.

In this case, the Picard-Fuchs system is the set of 9 equations
\[
\left(\begin{array}{ccccccccc}
E & 0 & a & 0 & c/2 & -1/5 & 0 & 0 & 0\\
0& E & 0 & a & 0 & c/2 & -1/5  & 0 & 0\\
0 & 0 & E & 0 & a & 0 & c/2 & -1/5  & 0\\
0 & 0 & 0 & E & 0 & a & 0 & c/2 & -1/5 \\
0 & 2a & 0 & 2c & -1 & 0 & 0 & 0 & 0\\
0 & 0 & 2a & 0 & 2c & -1 & 0 & 0 & 0\\
0 & 0 & 0 & 2a & 0 & 2c & -1  & 0 & 0\\
0 & 0 & 0 & 0 & 2a & 0 & 2c & -1  & 0\\
0 & 0 & 0 & 0 & 0 & 2a & 0 & 2c & -1 
\end{array}\right)
\left(\begin{array}{c} I_0 \\ I_2 \\ I_3 \\ I_4 \\ I_5 \\ I_6 \\ I_7 \\ I_8 \\ I_9 \end{array}\right)
=
\left(\begin{array}{c} \zeta_0 \\ \mu_1 \\ \zeta_2 \\ \zeta_3 \\ 0 \\ 2\zeta_0 \\ 4\zeta_1 \\ 6\zeta_2 \\ 8\zeta_3\end{array}\right),
\]
where
\[
\zeta_k:=\oint\frac{2 v^k~dv}{\sqrt{2(E+av^2+\frac{1}{2}cv^4-v^5)}},\qquad k=0,1,2,\ldots
\]
and $I_k$ are defined similarly; see \cite{Bronski.Johnson.Kapitula.2011}, for instance. 
In particular, notice that $T=\zeta_1$, $M=\zeta_3$ and $P=\zeta_5$.  
While the expressions are too long to list here, they are easily handled with the aid of a computer
algebra system, in which case one finds that $\Delta_{MI}>0$ in the entire domain of existence.  
In conclusion,
\emph{all periodic traveling waves of \eqref{kdv1/2} are modulationally stable}; see \cite{Bronski.Johnson.Kapitula.2011} for details.

\subsection{Extensions to Equations with Nonlocal Dispersion}

The analysis in previous sections readily extends to other equations allowing for \emph{nonlocal} dispersion.
We will illustrate this by discussing equations of KdV type 
\begin{equation}\label{nlkdv2}
u_t=\mathcal{M}u_x+f(u)_x,
\end{equation}
where $\mathcal{M}$ is defined via the Fourier series 
\[
u(x)=\sum_{n\in\mathbb{Z}}\hat{u}(n)e^{inx}\quad\Rightarrow\quad
(\mathcal{M}u)(x)=\sum_{n\in\mathbb{Z}}m(n)\hat{u}(n)e^{inx},
\]
characterizing dispersion in the limit limit. We assume throughout that the symbol $m$ is even and real-valued.

In the case of $m(k)=-k^2$, notably, \eqref{nlkdv2} recovers the generalized KdV equation.
Notice, however, that \eqref{nlkdv2} is nonlocal unless the dispersion symbol $m$ is a polynomial of $ik$. 
Examples include the KdV equation with fractional dispersion
\begin{equation}\label{fgkdv}
u_t=\Lambda^\alpha u_x+f(u)_x, \qquad -1\leq \alpha \geq 2,
\end{equation}
where $\Lambda=\sqrt{-\partial_x^2}$ is defined via its Fourier multiplier $m(n)=|n|$.
In the case of $\alpha=2$, it recovers the generalized KdV equation, and 
in the case of $\alpha=1$ and $f(u)=\frac12u^2$, it corresponds to the Benjamin-Ono equation,
which was proposed in \cite{Benjamin.1970} and \cite{Ono.1975}
as a model of the unidirectional propagation of internal waves of small amplitudes in deep water. 
Recalling that the Fourier multiplier of the Hilbert transform $\mathcal{H}$ is $-i\text{sgn}(k)$, 
the Benjamin-Ono equation may alternatively be written 
\begin{equation}\label{gbo}
u_t=\mathcal{H}\partial_x^2u+uu_x.
\end{equation}
Like the KdV and modified KdV equations, the Benjamin-Ono equation is completely integrable,
although, of course, our methods do not rely on this fact. 
In the case of $\alpha=-1/2$ and $f(u)=\frac12u^2$, moreover, 
\eqref{nlkdv2} was argued in \cite{Hur.breaking} to have relevance to 
water waves in two dimensions and in the infinite depths. 
Incidentally fractional powers of the Laplacian occurs in many applications, 
such as dislocation dynamics in crystals (see \cite{CDLFM}, for instance)
and financial mathematics (see \cite{CT}, for instance). 
Notice that the dispersion symbol of \eqref{fgkdv} is \emph{homogeneous} 
in the sense that $m(\lambda k)=a^\alpha m(k)$ for all $k\in\mathbb{R}$ and $a>0$,
which makes the analysis for waves of all amplitudes tractable
by the methods described in the previous sections. 

However, many applications require non-homogeneous symbols
for their dispersive operators.  The intermediate long-wave (ILW) equation
\begin{equation}\label{ilw}
u_t+u_x+(1/H)u_x-\mathcal{N}_Hu_x+\left(u^2\right)_x=0.
\end{equation}
is such an example, where here $H>0$ is a parameter and the Fourier multiplier $\mathcal{N}_H$ is defined via its symbol $m(k;H)=k\coth(kH)$.
The ILW equation is indeed of the form \eqref{nlkdv2}, where the symbol is given by $m(\xi)=1+1/H-\mathcal{N}_H$, and it arises as a model
to describe the interface between density stratified fluids, both with finite depth \cite{Joseph.1977}.

Another equation with non-homgeneous dispersion corresponds to the choices $m(k):=\sqrt{\frac{\tanh(k)}{k}}$ and $f(u)=u^2$, in which case
\eqref{nlkdv2} becomes the Whitham equation.  The Whitham equation was originally put forward by Whitham 
in \cite[p. 477]{Whitham.book} to model unidirectional propagation of surface water waves with small amplitudes, but not necessarily long wavelengths,
in a channel.  Notice that
\[
\sqrt{\frac{\tanh(k)}{k}}=1-\frac{1}{6}k^2+\mathcal{O}(|k|^4),~~|k|\ll 1
\]
so that we may regard the KdV equation \eqref{kdv} as to approximate up to ``second" order the dispersion relation for the Whitham equation,
and hence for the full water wave problem, in the long-wavelength regime (i.e. small frequencies $k$).  Unlike the KdV, however,
the Whitham equation may be used to model short wave uni-directional phenomena.  As a matter of fact, Whitham advocated
that this equation could be used to explain wave breaking and peaking, both of which are unobserved in the KdV approximation.
In light of the famous Benjamin-Feir instability, it is natural to expect that small periodic wave trains of the Whitham equation are modulationally
unstable provided their wavelength is sufficiently short.

Unlike the gKdV equation \eqref{gkdv}, the existence of periodic traveling waves for nonlocal equations of the form \eqref{nlkdv2}
is not trivial.  Indeed, a traveling wave of \eqref{nlkdv2} is a solution $u(x,t)=u(x-ct)$ where $u(\cdot)$ satisfies, after
one integration, the profile equation
\begin{equation}\label{nlquad}
-\mathcal{M}u+cu+f(u)=a
\end{equation}
for some constant $a\in\mathbb{R}$. 
For a broad range of dispersion symbols and nonlinearities, 
a plethora of periodic traveling waves of \eqref{nlkdv2} may be attained from variational arguments,
e.g. the mountain pass theorem applied to a suitable functional whose critical points satisfy \eqref{nlquad}.
In each example considered below, there will exist a smooth four-parameter (including translations) family of periodic traveling waves,
in agreement with the existence theory in the local gKdV case considered in previous sections.
%
%

Once existence and an appropriate parametrization of the periodic solutions of \eqref{nlquad} has been established,
modulational stability can be studied through studying the $L^2(\mathbb{R})$-spectrum of the linearized operator
\[
L:=\partial_x\left(-\mathcal{M}+c-f'(u)\right)
\]
near the origin.  As in the case of local gKdV equations, since the function $u$ is $T$-periodic, the $L^2(\mathbb{R})$-spectrum of $L$ is 
readily described by Floquet-Bloch theory:
\[
{\rm spec}_{L^2(\mathbb{R})}(L)=\bigcup_{\xi\in[-\pi/T,\pi/T)}{\rm spec}_{L^2_{\rm per}([0,T])}(L_\xi),
\]
the Bloch operators $L_\xi:=e^{-i\xi x}Le^{i\xi x}$ act as usual on $L^2_{\rm per}([0,T])$.
The modulational stability theory developed in Section \ref{s:theory} can thus be applied directly
in this case \emph{provided} the Jordan structure of the generalized kernel of $L$ is analogous to that described in Lemma \ref{l:gker}.

\subsubsection{The Benjamin-Ono Equation}\label{s:bo}

Let us consider the Benjamin-Ono equation (after the change of variables $u\mapsto -u$ in \eqref{gbo})
\begin{equation}\label{bo}
u_t-\Lambda u_x+(u^2)_x=0.
\end{equation}
In \cite{Benjamin.1970}, Benjamin exploited the Poisson summation formula 
and found that for each $k>0$ there exists a two-parameter (up to translations) family of $2\pi/k$-periodic traveling waves of \eqref{bo} of the form
\begin{equation}\label{bosoln}
u(z;a,k,c)=\frac{{\displaystyle \frac{k^2}{\sqrt{c^2-4a-k^2}}}}{{\displaystyle \sqrt{\frac{c^2-4a}{c^2-4a-k^2}}}-\cos(kz)}-\frac12(\sqrt{c^2-4a}+c),
\end{equation}
where $a$ and $c$ are arbitrary constants constrained by the conditions
\begin{equation}\label{boconstraint}
c<0\quad\textrm{and}\quad k^2<c^2-4a.
\end{equation}
It follows that \eqref{bo} admits a three parameter family (up to translations) of periodic traveling wave solutions, and that these solutions can be parametrized by the period and the parameters $a$ and $c$.  Here the period, equivalently $k$, is a parameter that is \emph{independent of $a$ and $c$} so that, in particular, derivatives of the solution
with respect to $a$ and $c$ are \emph{automatically $T$-periodic in $x$}.  
Since \eqref{bo} obeys Galilean invariance under
\begin{equation}\label{bogal}
u(z;a-c\lambda+s^2,k,c-2s) = u(z;a,k,c)+\lambda
\end{equation}
for all $\lambda\in\mathbb{R}$, upon an appropriate choice of $\lambda$,
one may assume that $a=0$. 

To analyze the modulational instability of a given solution $u$ of the above form,  we turn our attention to the 
$L^2_{\rm per}([0,T])$ spectrum near the origin of the Bloch operators 
\[
L_\xi:=e^{-i\xi x}\partial_x\left(\Lambda+c-2u\right)e^{i\xi x}
\]
for $|\xi|\ll 1$.  Although the symbol of $\Lambda$ is not smooth near the origin,
the symbol of $\partial_x\Lambda$ is $C^1$ at the origin and we can expand the Bloch operators
as
\[
L_\xi=L_0+(i\xi)L_1+(i\xi)^2 L_2+o(|\xi|^2)
\]
where 
\begin{equation}\label{bolin}
L_0:=\partial_x\left(\Lambda-2u-c\right)=:\partial_x\mathcal{L}
\end{equation}
and
\[
L_1=[L_0,x]=2\Lambda-2u-c,\qquad L_2=[L_1,x]=\Lambda^{-1}\partial_x.
\]
Following the analysis of Section \ref{s:theory}, the first step in modulational stability analysis is to understand the Jordan structure
of the kernel of the unmodulated operator $L_0$.  This is achieved in the following lemma.

\begin{lemma}[\cite{Bronski.Hur.2014}]\label{boker}
Suppose that $u=u(\cdot;a_0,k_0,c_0)$ is a $2\pi/k$-periodic traveling wave solution of \eqref{bo} and that
$\{M,P\}_{a,c}:=\det\left(\frac{\partial(M,P)}{\partial(a,c)}\right)$ is non-zero at $(a_0,k_0,c_0)$.  Then the 
generalized $L^2_{\rm per}([0,T])$ kernel of the linear operator $L_0$ defined in \eqref{bolin} possesses the following Jordan block structure:
\begin{itemize}
\item[(i)]~ ${\rm dim(ker}(L_0))=2$
\item[(ii)]~ ${\rm dim(ker}(L_0))/{\rm dim(ker}(L_0))=1$
\item[(iii)]~ ${\rm dim(ker}(L_0^{n+1}))/{\rm dim(ker}(L_0^n))=0$ for all integers $n\geq 2$.
\end{itemize}
In particular, the functions
\[\begin{array}{lcl}
v_0=u_a & \qquad & w_0=M_c u-P_c\\
v_1=u_x & \qquad & w_1=\int_0^x\left(M_au_c-M_cu_a\right)dx\\
v_2=u_c & \qquad & w_2=P_a-M_a u
\end{array}
\]
form a basis and dual basis for the generalized kernel of $L$, respectively:
\begin{align*}
&L_0v_0=Lv_1=0,\quad L_0v_2=v_2\\
&L_0^\dag w_0=L^\dag w_2=0,\quad L_0^\dag w_1=w_2.
\end{align*}
Moreover, $\left<w_j,v_k\right>=-\{M,P\}_{a,c}\delta_{jk}$ and $\left<w_j,L_lv_k\right>=0$ whenever $j+k+l$ is even.
\end{lemma}


The proof follows the same lines as the proof of Lemma \ref{l:gker}; see \cite{Bronski.Hur.2014}.  The key step to extend the proof to the current nonlocal case
is to verify that the kernel of $L_0$ is only two dimensional.  This, would follow essentially immediately provided one could verify that the operator
$\mathcal{L}$ in \eqref{bolin} is non-degenerate, i.e. that
\begin{equation}\label{nondeg}
{\rm ker}\left(\mathcal{L}\right)={\rm span}\{u_x\}.
\end{equation}
While this is obvious in the local gKdV case, by classical Sturm-Liouville theory, it is far from obvious in this more general nonlocal context.
Nevertheless, by relating via complex analysis techniques the nonlocal profile equation \eqref{nlquad} to a fully nonlinear ODE,
Amick and Tolland \cite{Amik.Toland.1991} demonstrated that the operator $\mathcal{L}$ associated to the Benjamin-Ono equation
is indeed non-degenerate in both the solitary and periodic settings.  Equipped with this non-degeneracy result, the 
proof of Lemma \ref{boker} follows from the same arguments as Lemma \ref{l:gker}.

To verify the hypothesis of Lemma \ref{boker}, owing to the explicit solution formulae \eqref{bosoln}, 
we note that a straightforward calculation (see \cite{Bronski.Hur.2014} for details) reveals that 
\begin{align*}
M(a,k,c):=&\int_0^{2\pi/k}u(z;a,k,c)~dz=2\pi-\frac{\pi}{k}(\sqrt{c^2-4a}+c)
\intertext{and}
P(a,k,c):=&\frac{1}{2}\int_0^{2\pi/k}u^2(z;a,k,c)~dz=-c\pi+\frac{\pi}{4k}(\sqrt{c^2-4a}+c)^2.
\end{align*}
so that, in particular,
\[
\{M,P\}_{a,c}=\frac{2\pi^2}{k\sqrt{c^2-4a}}>0
\]
for all $(a,k,c)$ satisfying the constraints \eqref{boconstraint}.
Lemma \ref{boker} thus provides an explicit basis for the generalized kernels of the unmodulated operator
$L_0$ and its adjoint about any periodic traveling wave solution of \eqref{bo}.

With Lemma \ref{boker} in hand, the methodology of Section \ref{s:migkdv} may be applied \emph{without modification}.
In particular, the triple eigenvalues at $\lambda=0$ of $L_0$ bifurcates into 3 eigenvalues $\lambda_j(\xi)$ that are $C^1$
in $\xi$ near $\xi=0$ and satisfy
\[
\lambda_j(\xi)=i\mu_j\xi+o(\xi),~~j=1,2,3
\]
where the constants $\mu_j\in\mathbb{C}$ are the eigenvalues of, taking into account the parametrization of solutions in \eqref{bosoln}, 
the effective dispersion matrix ${\bf D}(a,k,c)$ in \eqref{edm}.  
Using \eqref{bosoln} this matrix can be calculated \emph{explicitly} in terms of the traveling wave parameters $(a,k,c)$
\emph{provided} an explicit formula for the corrector $v_1^{(1)}$ in \eqref{basisexpand} can be found: due to the parametrization of
solutions of \eqref{bo} by their period, the formula from our gKdV analysis must be modified.  Since \eqref{Po} still
holds in the present case, $v_1^{(1)}$ may be identified as an appropriate multiple of the $2\pi/k$ periodic
function $xu_x+Tu_T$ where $T:=\frac{2\pi}{k}$; see \cite[Lemma 9]{Bronski.Hur.2014}.
With this choice, the matrix ${\bf D}(a,k,c)$ can be shown, 
up to a similarity transformation\footnote{Indeed, this formula for ${\bf D}$ was derived in \cite{Bronski.Hur.2014} using a direct spectral
perturbation expansion, a method equivalent to that described in Section \ref{s:migkdv} above.}
to be given (for $a=0$) by
\[
{\bf D}(a=0,k,c)=
\left(\begin{array}{ccc}
-\pi T & (\pi T)^2-(\pi/c)^2 & 0\\
1 & \pi T & 0\\
2\pi^2 & 0 & \pi T
\end{array}\right),\quad T=\frac{2\pi}{k}.
\]
A quick calculation shows the eigenvalues of ${\bf D}(0,k,c)$ are $\pi T$ and $\pm\pi T\sqrt{2-(cT)^{-2}}$.  Since
periodic traveling waves of \eqref{bo} exist only when $c^2>k=1/T$, it follows that the eigenvalues
of ${\bf D}(0,k,c)$ are real and distinct.  By the Galilean invariance \eqref{bogal} and Theorem \ref{t:whitham}, this yields
the following result.

\begin{theorem}[\cite{Bronski.Hur.2014}]
All periodic traveling wave solutions of the Benjamin-Ono equation \eqref{bo} are modulationally stable.
\end{theorem}

In fact, the above calculations on the Benjamin-Ono equation \eqref{bo} is merely a special case of the recent
analysis of Bronski and Hur where they considered the modulational instability of arbitrary amplitude periodic
traveling wave solutions to the fractional KdV equation \eqref{fgkdv} with $1/3<\alpha<2$.  For such $\alpha$,
one can show through the calculus of variations that for each $T>0$ there exists a two-parameter (up to translations) family of periodic
traveling wave solutions parametrized by the constant of integration $a$ and the wave speed $c$ and that, furthermore, these
waves depend on $a$, $T$, and $c$ in a $C^1$ manner.  A key difficulty in this more general setting is to verify the non-degeneracy
property \eqref{nondeg} for all $1/3<\alpha\leq 2$.   Although the methods of Amik \& Tolland \cite{Amik.Toland.1991}, described above, do not
apply outside the Benjamin-Ono case $\alpha=1$, it is possible to verify through a periodic adaptation of the 
recent nonlocal Sturm-Liouville theory of Frank \& Lenzmann \cite{Frank.Lenzmann.2013} for fractional Schrodinger operators to
verify that the kernel of $L_0$ is indeed two dimensional for all $1/3<\alpha\leq 2$; see \cite{Bronski.Hur.2014,Hur.Johnson.2015} for details.  Equipped
with this non-degeneracy result, the effective dispersion matrix can be numerically computed for all such $\alpha$; see\cite{Bronski.Hur.2014} for details.

\subsubsection{The Whitham Equation for Water Waves}\label{s:wsmall}

As described above, in the absence of explicit solution formulas it is in general very difficult to calculate the effective dispersion
matrix \eqref{edm} explicitly, and one must typically rely on well-conditioned numerical techniques.  One case
where this difficulty is avoided is in the small amplitude limit.  Indeed, as we will see below, for a wide class of nonlocal evolution
equations of the form \eqref{nlkdv2} such small amplitude solutions can be expanded (in amplitude) explicitly to arbitrarily high algebraic order.
Equipped with such an explicit expansion of solutions, one can compute an aymptotic expansion of the now amplitude dependent
effective dispersion matrix \eqref{edm} and hence ascertain the modulational stability for waves of sufficiently small amplitudes.

To illustrate this procedure, we now consider the modulational stability of periodic traveling wave solutions of the Whitham equation
\begin{equation}\label{Weqn}
u_t+\mathcal{M}u_x+\left(u^2\right)_x=0
\end{equation}
where here the symbol of the operator $\mathcal{M}$ is given by $m(n):=\sqrt{\frac{\tanh(n)}{n}}$.  As described above, equation \eqref{Weqn}
combines the full uni-directional dispersion relation for surface water waves with a canonical shallow water nonlinearity.
Details of the forthcoming analysis can be found in \cite{Hur.Johnson.2015b}.

Periodic traveling wave solutions of \eqref{Weqn} take the form $u(x,t)=u(x-ct)$, where $c>0$ and $u$ satisfies
\begin{equation}\label{Wquad}
\mathcal{M}u-cu+u^2=(1-c)^2b
\end{equation}
for some constant $b\in\mathbb{R}$; the normalizing factor of $(1-c)^2$ is added for convenience.  To describe the small amplitude periodic wave trains of \eqref{Weqn}
we seek solutions of \eqref{Wquad} of the form $u(x)=w(z)$ where $z=k x$, where here $w$ is required to be $2\pi$-periodic
and $k>0$ is interpreted as a wave number.  The $2\pi$-periodic function $w$ is then required to be a solution of
\begin{equation}\label{Wquad2}
\mathcal{M}_k w-cw+w^2=(1-c)^2b
\end{equation}
where here $\mathcal{M}_k$ is defined via
\[
\mathcal{M}_k e^{inz}=m(k n)e^{inz},~~n\in\mathbb{Z}.
\]
Observe that the symbol here is \emph{not homogeneous} and hence the wavenumber may not be factored out of the symbol.

Noting that $\mathcal{M}_k:H^s_{2\pi}\to H^{s+1/2}_{2\pi}$ is a continuous operator for all $k>0$ and $s\geq 0$, a simple
argument verifies that solutions of \eqref{Wquad2} with $\|w\|_{L^\infty}<\frac{c}{2}$ belong to the Sobolev space $H^\infty_{2\pi}$ and hence 
are smooth.  As a result, smooth $2\pi$-periodic solutions of \eqref{Wquad2} with small amplitudes may be sought as roots of the function
\[
F(w,k,c,b)=\mathcal{M}_k w-cw+w^2-(1-c)^2b;
\]
note that $F:H^1_{2\pi}\times\mathbb{R}_+\times\mathbb{R}_+\times\mathbb{R}\to H^1_{2\pi}$ by a Sobolev inequality.
Clearly $F(0;k,c,0)=0$ for every $k,c\in\mathbb{R}_+$ and, furthermore, the derivative of this map at such a trivial solution
is the operator
\[
\partial_w F(0,k,c,0)=\mathcal{M}_k -c
\]
which has a trivial kernel in $H^1_{2\pi}$ provided that $c\neq m(k n)$ for any $n\in\mathbb{Z}$.  By the implicit function theorem, 
near any such trivial solution for which $\partial_wF$ has a trivial kernel, there exists only other trivial solutions.
Thus, to construct non-trivial solutions of \eqref{Wquad2} we  start at a zero solution for which the kernel of $\partial_wF$ is non-trivial.

To this end, notice that if $c=m(k)$ then
\[
{\rm ker}\left(\partial_wF(0,k,c,0)\right)={\rm span}\left\{e^{\pm iz}\right\}.
\]
Using a straightforward Lyapunov-Schmidt reduction, it follows that for each $k>0$ there exists a two-parameter (up to translations) family of non-trivial roots of $F$,
corresponding to smooth (even) $2\pi$-periodic solutions of \eqref{Wquad2}.  Furthermore, their small amplitude asymptotics may be found by seeking
solutions of \eqref{Wquad2} of the form
\begin{align*}
w(k,A,b)(z)&=w_0(k,b)+A\cos(z)+A^2w_2(z,b)+A^3w_3(z,b)+\mathcal{O}(|A|^4)\\
c(k,A,b)&=m(k)+A^2c_2(k,b)+\mathcal{O}(|A|^4).
\end{align*}
where here $A$ is a real parameter (a sort of generalized amplitude), the $w_j$ are taken to be even, and $Q_b$ is defined for $|b|\ll 1$ via $F(w_0(k,b),k,m(k),b)=0$.
Plugging the above expansions into \eqref{Wquad2} we find
\begin{equation}\label{Wsolnexpand}
\begin{aligned}
w(k,A,b)(z)&=w_0(k,b)+A\cos(z)+\frac{1}{2}A^2\left(\frac{1}{m(k)-1)}+\frac{\cos(2z)}{m(k)-m(2k)}\right)\\
&\qquad+\mathcal{O}(A(A^2+b^2))\\
c(k,A,b)&=c_0(k,b)+A^2\left(\frac{1}{m(k)-1}+\frac{1}{2(m(k)-m(2k))}\right)\\
&\qquad+\mathcal{O}(A(A^2+b^2))
\end{aligned}
\end{equation}
as $|A|,|b|\to 0$, where
\[
c_0(k,b)=m(k)+2b(1-m(k))-6b^2(1-m(k))+\mathcal{O}(b^3)
\]
and
\[
w_0(k,b)=b(1-m(k))-b^2(1-m(k))+\mathcal{O}(b^3);
\]
although higher order expressions are obtainable, they are not necessary in our analysis. 
Together with translation invariance, this yields a four parameter family of small amplitude periodic traveling waves of the Whitham equation \eqref{Weqn}:
see \cite{Hur.Johnson.2015b,Ehrnstrom.Kalisch.2009,Ehrnstrom.Kalisch.2013} for more details.
Throughout the remainder, we take $b=0$ for simplicity.

With the above parametrization in hand, we now turn to study the $L^2_{2\pi}:=L^2_{\rm per}([0,2\pi])$ spectrum near the origin of the Bloch 
operators
\[
L_\xi(k,A):=e^{-i\xi z}L(k,A)e^{i\xi z},\quad L(k,A):=\partial_z\left(-\mathcal{M}_k+c(k,A,0)-2w(k,A,0)\right)
\]
for $|(\xi,A)|\ll 1$.
Unlike the case of the Benjamin-Ono equation considered in Section \ref{s:bo} above, the symbol $m(\cdot)$ is smooth near the origin
and hence we may easily expand the Bloch operators as
\[
L_\xi(k,A)=L_0(k,A)+(i\xi)L_1(k,A)+(i\xi)^2L_2(k,A)+\mathcal{O}(|\xi|^2)
\]
where
\begin{align*}
L_0(k,A)&=\partial_z\left(-\mathcal{M}_k+m(k)\right)-2a\partial_z\left(\cos(z)\cdot\right)+\mathcal{O}(A^2+b^2)\\
L_1(k,A)&=\left[-\partial_z\mathcal{M}_k,z\right]+m(k)-2a\cos(z)+\mathcal{O}(A^2+b^2)\\
L_2(k,A)&=\left[\left[-\partial_z\mathcal{M}_k,z\right],z\right]+\mathcal{O}(A^2+b^2)
\end{align*}
for $|A|\ll 1$.  
As in our previous analysis, the key step in the modulational stability analysis is to understand the Jordan structure of the kernel
of the unmodulated operator $L_0(k,A)$.  
To study the spectrum of $L_\xi$ near the origin, we first note that the operator estimate
\[
\left\|L_{\xi}(k,A)-L_\xi(k,0)\right\|_{L^2_{2\pi}\to L^2_{2\pi}}=\mathcal{O}(A)
\]
as $A\to 0$, uniformly in $\xi\in[-1/2,1/2)$, which may be verified by bruit force, implies that the spectrum of $L_\xi(k,A)$ is close to that of 
the constant coefficient operator $L_\xi(k,0)$, uniformly in the Bloch frequency $\xi$.  Consequently, we may consider the spectrum
of $L_\xi(k,A)$ as a small perturbation of that of $L_\xi(k,A)$.  In particular, the spectral projections $\Pi_\xi$ and
$\widetilde\Pi_\xi$ from Section \ref{s:theory} are now \emph{amplitude dependent}, described to leading order by
the projections associated to the trivial solution $A=0$.  

We concentrate first on the right eigenprojection.
Being constant coefficient, the spectrum of $L_\xi(k,0)$ for each $k>0$ and $\xi\in[-1/2,1/2)$ can be computed via Fourier analysis
as
\[
{\rm spec}\left(L_\xi(k,0)\right)=\left\{i\omega_{n,\xi}(k):n\in\mathbb{Z}\right\}\subset i\mathbb{R}
\]
where $\omega_{n,\xi}(k):=(n+\xi)\left(m(k)-m(k n+k\xi)\right)$ with corresponding eigenfunctions $e^{inz}$.
In particular, when $\xi=0$ we have $\omega_{-1,0}=\omega_{0,0}=\omega_{1,0}=0$ so that $\lambda=0$ is an eigenvalue of $L_0(k,0)$
of algebraic multiplicity three, with eigenfunctions $\cos(z)$, $\sin(z)$, and $1$.  Further, it follows from above that the eigenspaces
associated to the three non-zero eigenvalues $\omega_{\pm 1,\xi}, \omega_{0,\xi}$ are independent of $\xi$.  As such,
the associated right spectral projection $\Pi_\xi(k,0)$ is independent of $\xi$.  Finally, we note that the adjoint
operators $L_\xi(k,0)$, being constant coefficient, share the same eigenvalues and eigenfunctions, and hence 
the left spectral projection $\widetilde\Pi_\xi(k,0)$ is independent of $\xi$ as well.  In particular, we find\footnote{Following \cite{Hur.Johnson.2015b}, 
the projections here are slightly rescaled from that of Section \ref{s:theory}.  Here we use 
$M_\xi:=\left[\langle\psi_i,L_\xi\phi_j\rangle/\langle\psi_i,\phi_j\rangle_{L^2_{\rm per}([0,T])}\right]_{i,j=1,2,3}$
and similarly for the projection of the identity.}
\[
M_0(k,0):=\left(\Pi_\xi L_\xi\Pi_\xi\right)(k,0)=\left(\begin{array}{ccc}
															\frac{i}{2}(\omega_{1,\xi}+\omega_{-1,\xi}) & \frac{1}{2}(\omega_{1,\xi}-\omega_{-1,\xi}) & 0\\
															- \frac{1}{2}(\omega_{1,\xi}-\omega_{-1,\xi}) & \frac{i}{2}(\omega_{1,\xi}+\omega_{-1,\xi}) & 0\\
															0 & 0 & i\omega_{0,\xi}
															\end{array}\right).
\]

For $|A|\ll 1$, it follows that the spectral projections can be expanded as
\begin{align*}
\Pi_\xi(k,A)=\Pi_0(k,A)+\mathcal{O}(A^2+\xi^2),~~\widetilde\Pi_\xi(k,A)=\widetilde\Pi_0(k,A)+\mathcal{O}(A^2+\xi^2)
\end{align*}
so that, to first order, the spectral projections for the left and right eigenspaces of $L_\xi(k,A)$ agree with those of the unmodulated
operator $L_0(k,A)$.  In particular, the variations in these projections in $\xi$ do not effect the leading order asymptotics and, furthermore,
the bases used to construct these projections agree to first order.  As such, to construct the 
projections $\Pi_\xi(k,A)$ and $\widetilde\Pi_\xi(k,A)$ to leading order, we need only find a basis for the generalized kernel of $L_0(k,A)$
The origin is an eigenvalue of $L_0(k,A,b)$ of algebraic multiplicity
three and geometric multiplicity two for all $|A|,|b|\ll 1$ and, furthermore, the functions
\begin{align*}
\phi_1(z):=&\frac{1}{2(1-m(k))}\left[(\partial_b c)\partial_a w-(\partial_a c)\partial_b w\right](k,A,b)(z)\\
=&\cos(z)+\frac{-1/2+\cos(2z)}{m(k)-m(2k)}~A-6\cos(z)b+\mathcal{O}(A^2+b^2),\\
\phi_2(z):=&-\frac{1}{a}\partial_zw(k,A,b)(z)\\
=&\sin(z)+\frac{\sin(2z)}{m(k)-m(2k)}~A+\mathcal{O}(A^2+b^2)\\
\phi_3(z):=&1
\end{align*}
form a basis of the generalized eigenspace for $L_0(k,A,b)$ for all $|A|,|b|\ll 1$; specifically
\[
L_0\phi_1=L_0\phi_2=0,\quad\textrm{and}\quad L_0\phi_3=-2a\phi_2.
\]
Using this basis, we find
\begin{align*}
&M_\xi(k,A):=\left(\widetilde\Pi_\xi L_\xi\Pi_\xi\right)(k,A)=
\left(\begin{array}{ccc}
0 & 0 & 0\\ 0 & 0 & 2 \\ 0 & 0 & 0
\end{array}\right)
+i\xi\left(\begin{array}{ccc}
-k m'(k) & 0 & 0\\
0 & -k m'(k) & 0\\
0 & 0 & m(k)-1\end{array}\right)\\
&\quad-i\xi A\left(1+\frac{m(k)-1}{2(m(k)-m(2k))}\right)
\left(\begin{array}{ccc}
0 & 0 & 2\\ 0 & 0 & 0\\ 1 & 0 & 0
\end{array}\right)\\
&\qquad+\xi^2\left(k m'(k)+\frac{1}{2}k^2m''(k)\right)
\left(\begin{array}{ccc}
0 & -1 & 0\\ 1 & 0 & 0\\ 0 & 0 & 0\end{array}\right)+\mathcal{O}(|\xi|^3+|A|^3)
\end{align*}
and, similarly, that
\[
\left(\widetilde\Pi_\xi\Pi_\xi\right)(k,A)={\bf I}
-\frac{A}{m(k)-m(2k)}
\left(\begin{array}{ccc}
0 & 0 & 1\\ 0 & 0 & 0\\ 1/2 & 0 & 0
\end{array}\right)+\mathcal{O}(A^2)
\]

Next, we turn to study the roots of the characteristic polynomial
\begin{align*}
P(\lambda,\xi;k,A)&=\det\left(M_\xi(k,A)-\lambda\left(\widetilde\Pi_\xi\Pi_\xi\right)(k,A)\right)\\
&=c_3(\xi;k,A)\lambda^3+ic_2(\xi;k,A)\lambda^2+c_1(\xi;k,A)\lambda+ic_0(\xi;k,A)
\end{align*}
for each fixed $k>0$ and $(\lambda,\xi,A)|\ll 1$.  As before, the roots of $P(\cdot,\xi;k,A)$ correspond to the eigenvalues
of the Bloch operators $L_\xi(k,A)$ bifurcating from the origin for $|(\xi,A)|\ll 1$.  We first note that the coefficient
functions $c_j$ depend smoothly on $\xi$ and $A$ for $|(\xi,A)|\ll 1$ and that, moreover, the $c_j$ are real-valued since
the spectrum of $L_\xi(k,A)$ is symmetric about the imaginary axis.  Further, the spectral symmetry
\[
{\rm spec}_{L^2_{2\pi}}(L_\xi(k,A))=\overline{{\rm spec}_{L^2_{2\pi}}(L_{-\xi}(k,A))}
\]
implies that the functions $c_3$ and $c_1$ are even functions of $\xi$ while $c_2$ and $c_0$ are odd, and that the $c_j$'s are all even in $A$.
Since $\lambda=0$ is clearly a root of $P(\lambda,0;k,A)$ with multiplicity three and $\xi=0$ is a root of $P(0,\xi;k,A)$ with multiplicity
three for all $|A|\ll 1$, it follows that 
\[
c_j(\xi;k,A)=d_j(\xi;k,A)\xi^{3-j},\quad j=0,1,2,3
\]
for some real functions $d_j$ that depend smoothly on $A$ and $\xi$ and are even in $A$.  Thus, the roots of $P(\cdot,\xi)$ may be 
written as $\lambda=-i\xi X$ where $X$ is a root of 
\begin{align*}
&(i\xi)^{-3}P(-i\xi X,\xi;k,A)=\det\left(\frac{1}{i\xi}M_\xi(k,A)+X\left(\widetilde\Pi_\xi\Pi_\xi\right)(k,A)\right)\\
&\qquad\qquad=-d_3(\xi;k,A)X^3+d_2(\xi;k,A)X^2+d_1(\xi;k,A)X-d_0(\xi;k,A);
\end{align*}
compare to \eqref{charinform}.
It follows  that the underlying periodic traveling wave is modulationally unstable if the characteristic polynomial $P(-i\xi\cdot,\xi;k,A)$
admits a pair of complex roots, i.e. if its discriminant
\[
\Delta_{\xi,k,A}:=\left(18d_3d_2d_1d_0+d_2^2d_1^2+4d_2^3d_0+4d_3d_1^3-27d_3^2d_0^2\right)(\xi,k,A)<0
\]
for $|(\xi,A)|\ll 1$, and that it is modulationally stable if the polynomial admits three real roots, i.e. if $\Delta_{\xi,k,A}>0$.

Since the $d_j$'s are even in $A$, we may expand
\[
\Delta_{\xi,k,A}=\Delta_{\xi,k,0}+\Lambda(k)A^2+\mathcal{O}(A^2(A^2+\xi^2))
\]
as $|(\xi,A)|\ll 1$, expressing that, to leading order, the modulational stability is governed by that of the limiting constant state at $A=0$.
Using the above explicit formulae, we find that $\Delta_{0,k,0}=0$ and
\[
\Delta_{\xi,k,0}=\frac{(\omega_{0,\xi}-\omega_{1,\xi})^2(\omega_{0,\xi}-\omega_{-1,\xi})^2(\omega_{1,\xi}-\omega_{-1,\xi})^2}{\xi^6}>0
\]
for all $|\xi|\ll 1$.  Therefore, the sign of the discriminant $\Delta_{\xi,k,A}$ for $|(\xi,A)|\ll 1$ is determined by
the sign of the function $\Lambda(k)$.  Indeed, from above we see that if $\Lambda(k)>0$ then $\Delta_{\xi,k,A}>0$
fro all $|(\xi,A)|\ll 1$, implying modulational stability, while if $\Lambda(k)<0$ then $\Delta_{\xi,k,A}<0$ for $|A|\ll 1$ fixed
and $|\xi|$ sufficiently small.  The modulational stability of the underlying small amplitude waves is thus determined completely
by the wavenumber $k$.

\begin{figure}
\begin{center}
\includegraphics[scale=0.8]{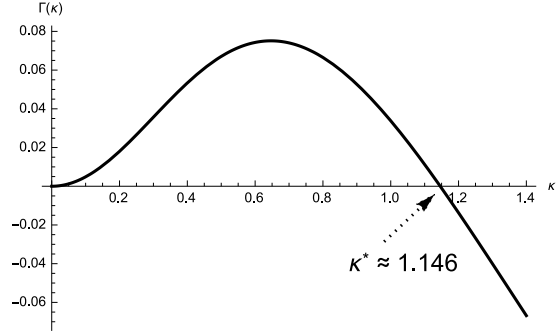}
\end{center}
\caption{Plot of $\Gamma(k)$ for Whitham equation.
}\label{WhithamFig}
\end{figure}

With the aid of computer algebra software, we find
\[
\Lambda(k)=\frac{(k(m(k)-1)'(k(m(k)-1)''}{m(k)-m(2k)}~\Gamma(k)
\]
where 
\[
\Gamma(k):=2(m(k)-m(2k))+(k(m(k)-1)'.
\]
Recalling that $m(k):=\sqrt{\frac{\tanh(k)}{k}}$, we see that the function $k\mapsto k(m(k)-1)$ is strictly
decreasing and concave down on the interval $(0,\infty)$, so that the sign of $\Lambda(k)$ is determined by the sign
of $\Gamma(k)$.  Further, a numerical evaluation of $\Gamma$ shows that there exists a unique $k^*\approx 1.146$
such that $\Gamma(k)>0$ for $k\in(0,k^*)$ and $\Gamma(k)<0$ for $k>k^*$; see Figure \ref{WhithamFig}.  

\begin{theorem}[\cite{Hur.Johnson.2015b}]\label{t:wmi}
A $2\pi/k$-periodic traveling wave solution of \eqref{Weqn} is sufficiently small amplitude is modulationally unstable
if $k>k^*\approx 1.146$, and modulationally stable if $0<k<k^*$.
\end{theorem}

Theorem~\ref{t:wmi} qualatatively captures the famous Benjamin-Feir instability in the context of the Whitham equation: 
see \cite{Hur.Johnson.2015b} and \cite{Stanford.et.al.2014} for more discussion on this connection.

\subsubsection{Fractional KdV and ILW Equations}

Finally, we note that the small amplitude analysis in Section \ref{s:whitham} is largely independent of the specific symbol for the nonlocal operator
$\mathcal{M}$.  Indeed, the only features of the symbol $m(\cdot)$ that was used in the derivation of the modulational stability
index $\Lambda(k)$ were (i) that $m(0)=1$ and (ii) that $k>0$ is such that $m(k)\neq m(nk)$ for all $n=2,3,4,\ldots$;
note that while (ii) provides a restriction on $k$, the condition $m(0)=1$ can be relaxed via a simple scaling argument\footnote{Although, the formula
for $\Lambda(k)$ must be modified accordingly.}.

Examples where this analysis can be applied \emph{without any modification} include that fractional KdV (fKdV) equation,
corresponding to \eqref{fgkdv} with $f(u)=u^2$ and the intermediate long-wave equation \eqref{ilw}.  For the fKdV, we find for each\footnote{The
condition $\alpha>1/2$ is an artefact of the corresponding existence theory.  See \cite{Johnson.2013} for details.} $\alpha>\frac{1}{2}$
that
\[
\Lambda_{\rm fKdV}(k;\alpha)=\frac{2k^{4\alpha}\alpha(1+\alpha)^4(2^{\alpha+1}-3-\alpha)}{2^{\alpha}-1},
\]
which is negative (positive) for all $k>0$ provided that $\alpha<1$ ($\alpha>1$), indicating modulational instability (stability) of
the associated small amplitude periodic wave trains. Recall that $\alpha=1$ corresponds to the Benjamin-Ono equation, which was
analyzed in Section \ref{s:bo}.  The corresponding result for more general power-law nonlinearities may be found in \cite{Johnson.2013}.

Similarly, for the ILW equation \eqref{ilw} we find for each $H>0$ that
\[
\Delta_{\rm ILW}(k;H)=\frac{(4H^2k^2-1)\cosh(Hk)+\cosh(3Hk)-8Hk\sinh(Hk))^2}{32H^4\sinh(Hk)^{12}}\Gamma_{ILW}(Hk)
\]
where $\Gamma_{ILW}(z):=1-2z^2-\cosh(2z)+2z\sinh(2z)$; see \cite{Hur.Johnson.2015b}.  
A numerical evaluation indicates that $\Gamma_{ILW}(z)>0$ for all $z>0$, and hence
that, for each $H>0$, the small amplitude periodic traveling wave solutions of \eqref{ilw} are modulationally stable\footnote{This rigorously justifies
a formal "amplitude equation" calculation in \cite{Pelinovsky.1995}.}.

\section{Concluding Remarks}\label{s:remarks}

We have made a detailed survey of many recent techniques rigorously studying 
the modulational instability of periodic traveling waves of KdV type equations. 
They are very general and apply to many other equations including systems involving dissipation.
They are independent of the integrability of the governing equations.  
The theoretical approach in Section~\ref{s:theory} may be used 
to rigorously validate formal predictions from Whitham's theory of modulations in a variety of settings; 
see \cite{Gavage.Noble.Rodrigues.2014}, \cite{Noble.Rodrigues.2013}, for instance,
for more concrete examples, other than those listed in Section~\ref{s:validate}.
We applied our theory to a number of concrete examples in Section~\ref{s:app}. 
To illustrate the robustness of our theory, we discussed extensions to nonlocal equations of KdV type 
and obtained detailed results for a number of concrete examples.  

We have been concerned with the linear (spectral) stability 
of periodic traveling waves of KdV type equations to slow modulations.
A fundamental open problem is to understand the \emph{nonlinear} dynamics 
associated with modulationally stable waves and, 
even more interestingly, induced by a modulational instability.  
They have been studied in the context of dissipative systems 
(see \cite{Johnson.Noble.Rodrigues.Zumbrun.2014}, for instance, and references therein), 
but it is not clear how such techniques apply to dispersive equations 
such as the gKdV equations. This is an exciting direction for future work.

\medskip
{\bf Acknowledgement}:
JCB is supported by the National Science Foundation grant DMS-1211364.
VMH is supported by the National Science Foundation grant CAREER DMS-1352597
and an Alfred P. Sloan Foundation Fellowship.
MAJ is supported by the National Science Foundation grant DMS-1211183.

\end{document}